\documentstyle[amstex,amssymb,11pt]{article}
  \input xypic
 \newcommand{\lon}{\longrightarrow}
 
 \newcommand{\rar}{\rightarrow}
 
 \newcommand{\Proof}{{\bf Proof}.\, }

 \newcommand{\Z}{{\Bbb Z}}
 \newcommand{\bS}{{\Bbb S}}
 \newcommand{\p}{{\partial}}

 \newcommand{\R}{{\Bbb R}}
 \newcommand{\ot}{\otimes}

 \newcommand{\Id}{\mbox{Id}}
 \newcommand{\Beq}{\begin{equation}}
 \newcommand{\Eeq}{\end{equation}}
 \newcommand{\Beqr}{\begin{eqnarray}}
 \newcommand{\Eeqr}{\end{eqnarray}}
 \newcommand{\Beqrn}{\begin{eqnarray*}}
 \newcommand{\Eeqrn}{\end{eqnarray*}}
 \newcommand{\Ba}{\begin{array}}
 \newcommand{\Ea}{\end{array}}
 \newcommand{\Bi}{\begin{itemize}}
 \newcommand{\Ei}{\end{itemize}}
 \newcommand{\Bc}{\begin{center}}
 \newcommand{\Ec}{\end{center}}
 \newcommand{\fg}{{\frak g}}

 \newcommand{\f}{{\cal O}}

 \newcommand{\caD}{{\cal D}}
 \newcommand{\cE}{{\cal E}}

 \newcommand{\cI}{{\cal I}}

 \newcommand{\cM}{{\cal M}}
 \newcommand{\cN}{{\cal N}}
 \newcommand{\cP}{{\cal P}}

 \newcommand{\al}{\alpha}
 \newcommand{\be}{\beta}
 \newcommand{\ga}{\gamma}

 \newcommand{\om}{\omega}

 \newcommand{\Img}{{\mathsf Im}\, }

 \newcommand{\sip}{\smallskip}
 \newcommand{\bip}{\bigskip}

\begin{document}

 \sloppy

 \title{Nijenhuis infinity and contractible dg manifolds\footnote{This work was
 partially supported by the G\"oran Gustafsson foundation.}}
 \author{ S.A.\ Merkulov}
 \date{}
 \maketitle
\begin{abstract}
We find a minimal differential graded (dg) operad  whose generic
representations in $\R^n$ are in one-to-one correspondence with
formal germs of those endomorphisms of the tangent bundle to $\R^n$
which satisfy the Nijenhuis integrability condition. This operad is of a surprisingly simple
origin --- it is the cobar
construction on the quadratic operad of homologically trivial dg Lie algebras.
As a by
product we obtain a strong homotopy generalization of this
geometric structure and show its homotopy equivalence to the structure of
contractible dg manifold.
\end{abstract}

\bip

\bip

{\bf 1. Introduction.} Extended deformation theory of a
mathematical structure $\frak X$ begins with the construction of
a ``controlling" algebra, $\fg$, over an operad $\caD ef$. The formal smooth
extended moduli space of deformations of
$\frak X$
 is then a homotopy class of minimal algebras over its
cofibrant resolution,
  $\caD ef_\infty$, canonically induced on the cohomology of $\fg$.
The operad $\caD ef$ contains often (if not always)
 the suboperad of Lie algebras which in
turn explains why a {part} of the data induced on the extended
moduli space can  also be effectively described with the help of the
classical idea of deformation functor on dg Artinian rings.
For example, if $\frak X$
is either the structure of associative algebra on a vector space
$V$, or a symplectic/complex structure on a manifold $M$, then the
appropriate $\caD ef_\infty$ is the minimal resolution of the
operad of Gerstenhaber algebras. A fact from which we can, for
example, immediately infer that the associated extended moduli
space (modelled as a formal supermanifold on the cohomology spaces
$Hoch^\bullet(A,A)$, $H^\bullet_{DRham}(M)$, and, respectively,
$H^\bullet_{Doulbeout}(M,\wedge^\bullet T_M)$, in the three cases just
mentioned) has canonically induced a Hertling-Manin$_\infty$
geometric structure \cite{HM,Me1}.

\sip

This machinery does not explain, however, what happens to the
structure $\frak X$ itself upon {\em generic}\, extended
deformation, that is, upon deformation in a direction other than  in
$Hoch^2(A,A)$, $H^2_{DRham}(M)$ or $H^1_{Doulbeout}(M,T_M)$
respectively. Thanks to Stasheff \cite{St} we know what
happens to associative algebra structure under generic
deformations--- it gets transformed into
$A_\infty$-structure. Remarkably enough, the latter can also be
understood in terms of minimal resolutions of operads!
This fact encourages our present attempt to understand some classical geometric
structures from  the operadic perspective. The hope is that the theory
of (di)operads can provide us with
an insight into what might be called {\em geometry}$_\infty$, the typical fibre
of  yet unknown ``universal curve" over generic point of the
extended moduli space.

\sip


In \cite{Me2} we presented a  dioperadic profile of
Poisson
geometry. In this paper we attempt to create  an operadic
profile of one of the key equations defining  complex structure, the
Nijenhuis integrability condition. After some preparations in
Sections 2-4, we show in the final Section 5 a
surprisingly simple minimal dg operad $N_\infty$ whose
generic representations in $\R^n$ are in one-to-one correspondence
with formal germs of those endomorphisms of the tangent
bundle, $J:T_{\R^n}\rar T_{\R^n}$, which satisfy the Nijenhuis
integrability condition, $N_J=0$. Generic representations of
$N_\infty$ in {\em graded}\, vector spaces provide us
then with the strong homotopy
generalization of this geometric structure. A minimal Nijenhuis$_\infty$ manifold
turns out to be a contractible dg manifold
in the sense of Kontsevich \cite{K}.

\bip

{\bf 2. Strongly homotopy Nijenhuis?} Let $M$ be a smooth manifold
and $J: T_M\rar T_M$ an endomorphism of the tangent bundle.
Nijenhuis \cite{N} used these data to construct a {tensor},
$$
\Ba{rccc} N_J: & \wedge^2 T_M & \lon & T_M \\
               & X\ot Y       & \lon & N_J(X,Y)
\Ea
$$
where
$$
N_J(X,Y):=[JX,JY] + J^2[X,Y] - J[X,JY] - J[JX,Y].
$$
By Newlander-Nirenberg theorem, its vanishing for an almost
complex structure $J$ is a necessary and sufficient condition for the latter
to be a complex structure. If $J:T_M\rar T_M$ is a fiberwise diagonalizable
endomorphism with all eigenvalues real and of constant multiplicity, then
eigenspaces
of $J$ are integrable if and only if $N_J=0$ \cite{Mi}.

\sip

A {\em Nijenhuis geometric structure}\, on $M$ is defined to be
an endomorphism $J: T_M\rar T_M$ satisfying the condition $N_J=0$.
Complex structures are, probably, the most important
examples of Nijenhuis  structures.

\sip

At present we have a  good experience  in obtaining
strongly homotopy generalizations of classical {\em algebraic}\, structures:
the major working tool is
the construction of  minimal dg
resolutions of the operads  governing that classical structures (see, e.g.
 \cite{GJ,GK,MSS}).
Can we apply this technique to get, say, a
 Nijenhuis$_\infty$ {\em geometric}\, structure? For that purpose let us
 first try to translate
 geometry into algebra, and
choose a local coordinate system, $\{t^\al\}_{1\leq \al \leq \dim M}$,
 at a point $*$
in $M$, that is a local isomorphism of germs of smooth
manifolds $(V=\R^n,0)\rar (M,*)$. The germ  at $*$ of any analytic endomorphism
$J: T_M\rar T_M$ gets represented explicitly as
\Beqrn
J &=& \sum_{\al,\be} J_\be^\al(t) dt^\be \ot \frac{\p}{\p t^\al},\\
J_\be^\al(t) &=&
\sum_{k\geq 0}\sum_{\ga_1,\ldots,\ga_k} J^\al_{\be\ga_1\ldots \ga_k}
t^{\ga_1}\ldots t^{\ga_k},
\Eeqrn
i.e.\ as an infinite collection of ``algebraic operations",
$\{J_k: \ot^{k+1}V\rar V\}_{k\geq 0}$. The Nijenhuis integrability constraint
would then translate into an infinite sequence of quadratic relations between
these operations. The experience with algebra, if understood straightforwardly,
suggests the following scenario: (1) consider an operad $\cN$ which governs
these algebraic operations and relations; (2) attempt to construct its
cofibrant resolution $\cN_\infty$;  (3) study generic
representations of the latter in dg vector spaces; in the particular
cases of vector spaces
concentrated in degree zero we are guaranteed to get all the operation $J_k$
satisfying the
required Nijenhuis integrability constraint (plus something else).

\sip

However, this scenario is  a doomed failure: the input operad $\cN$ is already so
monstrously
big that it is technically impossible (and repelling) to go for something yet bigger.

\sip

So let us try instead to go for something much smaller, and
investigate
a theoretical possibility that the Nijenhuis integrability condition
is already ``minimally resolved" in the sense that there exists
a minimal dg operad $N_\infty$
whose representations in spaces concentrated in degree zero are precisely
the collections $\{J_k: \ot^{k+1}V\rar V\}_{k\geq 0}$ satisfying the
required quadratic
constraints. Instead of steps (1) and (2) above we face now a problem of
reconstructing the full minimal dg operad $N_\infty$
from just its degree zero bit corresponding to $\cN$.

\sip

If the guess in the preceding paragraph is correct, then some operation(s) in
the family $\{J_k: \ot^{k+1}V\rar V\}_{k\geq 0}$ is (are) fundamental while all the
others are just higher homotopies of the fundamental one(s).
As the Nijenhuis integrability condition
is a system of {\em first}\, order non-linear partial differential
equations we suggest
the linear in $t^a$ bit, $J_1:\ot^{2}V\rar V$, of the endomorphism tensor $J$ for the role
of the fundamental operation, and immediately observe the following

\bip

{\bf 2.1. Lemma.} {\em There is a one-to-one correspondence,}
$$
\left\{\Ba{c}
{\mathrm Linear\ in}\ t^a \ \mathrm{ Nijenhuis\ structures}\\
             J: T_V\rar T_V
             \Ea\right\}
             \leftrightarrow
\left\{\Ba{c}
\mbox{pre-Lie algebra}\\
   {\mathrm structures\ on}\ V
             \Ea\right\}
$$

\bip

Thus pre-Lie algebras are going to play an important role in the story.

\bip

{\bf 2.2. Reminder on pre-Lie algebras.} A {\em pre-Lie}\, algebra is a vector space
together with a binary operation, $\circ: V^{\ot 2}\rar V$, satisfying the condition
$$
(a\circ b)\circ c - a\circ (b\circ c) - (-1)^{|b||c|}(a\circ c)\circ b +
(-1)^{|b||c|} a\circ (c\circ b) =0
$$
for any $a,b,c\in V$.

\sip

A pre-Lie algebra $V$ is naturally a Lie algebra with the bracket,
$$
[a,b]:=a\circ b  - (-1)^{|a||b|}b\circ a.
$$

The operad of pre-Lie algebras is Koszul \cite{CL} and its minimal dg resolution
can  be explicitly described using the technique of Ginzburg-Kapranov \cite{GK}.

\bip

{\bf 2.3. Remark.}
Pre-Lie$_\infty$ algebra structures will eventually form a (small) part of the
Nijenhuis$_\infty$ structure.
It turns out that to get the full picture of  Nijenhuis$_\infty$ one has to add to
the pre-Lie
operation $\circ$ at least one more algebraic operation, $[\, \bullet\, ]$,
 the Lie bracket operation
in degree 1.
Such an extension is more than expected: along the way from geometry to algebra
in Section 1 we fixed a flat structure at $(M,*)$, an ugly thing to do from the
geometric viewpoint. If everything is done properly, the resulting  Nijenhuis$_\infty$
geometric structure should not depend on such a choice, i.e.\ it must be tensorial.
Thus we must add an ingredient to our algebraic input whose homotopy theory
contains the group of general coordinate transformations as homotopy equivalences.
And Lie algebras in degree 1
just have this property: a minimal Lie$_\infty$ algebra structure is the same
as a homological vector field on a smooth formal manifold \cite{K};
generic diffeomorphisms
of the latter generate homotopy equivalence among the former.
Thus Lie algebra structure in degree 1 is expected to be present in the (di)operadic
profile of {\em any}\, diffeomorphism covariant geometric structure, cf.\
\cite{Me2}.

\bip

{\bf 3. Operad of pre-Lie$^{\mathbf\frak 2}$ algebras.}  Here is a required
graded extension of the notion
of pre-Lie algebra.

\sip

{\bf 3.1. Definition.} A {\em pre-Lie}$^{\mathbf\frak 2}$ algebra is a graded vector
space $V$ together with
two binary operations,
$$
\Ba{rccc}
\circ: &  V\ot V & \lon & V \\
       & a\ot b    & \lon & a\circ b
\Ea
\ \ \ \ \  , \ \ \ \
\Ba{rccc}
[\, \bullet\, ]: & \odot^2 V& \lon & V[1] \\
       & a\ot b    & \lon & (-1)^{|b|}[a\bullet b]
\Ea
$$
such that
\Bi
\item[(i)] the data $(V,\circ)$ makes $V$ into a pre-Lie algebra;
\item[(ii)] the data $(V, [\, \bullet\, ])$ makes $V[-1]$ into a Lie algebra; equivalently,
the odd Jacobi identity,
\Beqr\label{1}
[[a\bullet b]\bullet c]=[a\bullet[b\bullet c]] +
(-1)^{|b||a|+|b|+|a|}[b\bullet[a\bullet c]],
\Eeqr
holds for any $a,b,c\in V$;
\item[(iii)] a ``compatibility" identity,
$$
[a\bullet b]\circ c + (-1)^{|b|}a\circ [b\bullet c] + (-1)^{|b||a|+|b|}
b\circ [a\bullet c] = \hspace{7cm}
$$
\Beqr\label{2}
 \hspace{4cm}
 (-1)^{|b||c|+|c|} [(a\circ c)\bullet b]
+ (-1)^{(|a|+1)(|b|+|c|)+|a|} [(b\circ c)\bullet a].
\Eeqr
holds for any $a,b,c\in V$;
\Ei

\sip

The last condition (iii) should be explained.

\bip

{\bf 3.1.1. Interpretation I.} Recall that if
$(\fg, [\ ,\ ])$ is a Lie algebra and $M$
is a right Lie module over $(\fg, [\ ,\ ])$, then the Chevalley-Eilenberg construction
makes
$M\ot \wedge^\bullet \fg$ into a complex (computing $H_\bullet(\fg,M)$).

\sip

Similarly, if $(\fg, [\ \bullet \ ])$ is an odd Lie algebra (i.e. $[\, \bullet\, ]$
satisfies the odd Jacobi identity as in  (ii) above), and $M$ is  a module over
$(\fg, [\ \bullet \ ])$, then  the same Chevalley-Eilenberg construction makes
now the vector space $\odot^\bullet\fg\ot M$ into a complex.

\sip

Let $V$ be a graded vector space equipped with two operations $\circ$ and
$[\, \bullet\, ]$ satisfying conditions (i) and (ii) in definition 3.1. Then
the vector space $\odot^\bullet V \ot \wedge^\bullet V$ has two Chevalley-Eilenberg
 differentials:
\Bi
\item[(a)]
one Chevalley-Eilenberg differential, $d_\circ$, comes from the even Lie
algebra structure on $V$ (as
in Sect.\ 2.2)
and the right $(V,[\, ,\, ])$-module structure on $\odot^\bullet V$ generated by right
$\circ$-multiplication;
\item[(b)]
another Chevalley-Eilenberg differential, $d_\bullet$, comes from the odd Lie
algebra structure, $[\, \bullet\, ]$, on $V$
and the  $(V,[\, \bullet\, ])$-module structure on $\wedge^\bullet V$ generated by
the adjoint
$[\, \bullet\, ]$-action.
\Ei
The condition (iii) in the Definition 3.1 is then equivalent to commutativity of the
two  Chevalley-Eilenberg differentials, $d_\circ$ and $d_\bullet$.

\bip

{\bf 3.1.2. Interpretation II.} If $V$ is a vector space, then $V\oplus V[-1]$
is naturally a complex with trivial cohomology. If we write elements
of $V\oplus V[-1]$ as $a + \Pi b$, where $a,b\in V$ and $\Pi$ is a formal
symbol of degree 1,
 then the differential is given by
$$
d(a + \Pi b)= 0 + \Pi a.
$$

\sip

It is straightforward to check that a pair of linear maps, $\circ: \ot^2V\rar V$ and
$[\, \bullet\, ]: \odot^2 V\rar V[1]$,
 makes $V$ into a pre-Lie$^{\mathbf\frak 2}$-algebra
if and only if the brackets on $V\oplus V[-1]$,
\Beqrn
[a,b]&:=& a\circ b - (-1)^{|a||b|} b\circ a \\
\left[\Pi a, b\right] &:=& - (-1)^{|a|}[a\bullet b] +\Pi a\circ b\\
\left[\Pi a, \Pi b\right] &:=& \Pi\left[a\bullet b\right]
\Eeqrn
make the complex $(V\oplus V[-1], d)$ into a dg Lie algebra.

\bip

{\bf 3.1.3. Example.} Let $W$ be a vector space and $Free_{Lie}(W\oplus W[-1])$
the free Lie algebra on the vector space  $W\oplus W[-1]$. Consider two differentials,
$$
\Ba{rccc}
d: &  W\oplus W[-1] & \lon &  W\oplus W[-1]\\
       & a + \Pi b    & \lon & 0 +  \Pi a
\Ea
\ \   , \ \
\Ba{rccc}
q: &  W\oplus W[-1] & \lon &  W\oplus W[-1]\\
       & a + \Pi b    & \lon & b +  \Pi 0,
\Ea
$$
and extend them to $Free_{Lie}(W\oplus W[-1])$ as derivations of the Lie bracket.
For any Lie word $A$ of length $\lambda_{A}$ one has $(dq + qd)A=\lambda_A A$
so that after rescaling
$q \rar Q:= (\lambda_(\ldots))^{-1} q$ we get eventually two differentials on
$Free_{Lie}(W\oplus W[-1])$ which satisfy
$$
dQ + Qd =\Id.
$$
Hence $Free_{Lie}(W\oplus W[-1])$ splits canonically into the direct sum, $V\oplus V[-1]$,
where $V=\Img Q$. Using 3.1.2 it is now easy to check that
the  vector space $V$ is naturally
a
pre-Lie$^{\mathbf\frak 2}$ algebra with
\Beqrn
Qa \circ Qb&:=& Q[dQa, Qb] \\
\left[Qa\bullet  Qb\right] &:=& Q[dQa, dQb].
\Eeqrn
for any $a,b,c\in Free_{Lie}(W\oplus W[-1])$.
For example, if $w_1,w_2,w_3\in W\subset V$, then one has
\Beqrn
w_1 \circ w_2&=& \frac{1}{2}[w_1, w_2], \\
\left[w_1\bullet  w_2\right] &=&  \frac{1}{2}[w_1, \Pi w_2] - (-1)^{|w_1|}
 \frac{1}{2}[\Pi w_1, w_2],\\
 w_1\circ (w_2\circ w_3) &=&  \frac{1}{6}[w_1, [w_2, w_3]], \\
(w_1\circ w_2) \circ w_3 &=&  \frac{1}{3}[[w_1, w_2], w_3 ],
\Eeqrn
and so on.

\sip

In fact any Lie algebra with two derivations, $d$ and $Q$, of the Lie bracket
satisfying $d^2=Q^2=0$ and $dQ + Qd =\Id$ gives rise to a pre-Lie$^{\mathbf\frak 2}$
algebra structure
on $\Img Q\subset \fg$.

\bip

{\bf 3.1.4. Example.} Let $M$ be a smooth manifold and $\Omega^\bullet_M$ the associated
  graded commutative de Rham algebra of differentials forms. We shall show that derivations
 of the latter give naturally rise to one more example of pre-Lie$^{\mathbf\frak 2}$ algebra.

 \sip

 Instead of $M$ it is more
  suitable to introduce a $(\dim M|\dim M)$-dimensional
   supermanifold $\cM:=T_M[1]$, the total
 space of the tangent bundle to $M$ with the degree of typical fibre shifted by $1$.
 Then the de Rham algebra  $\Omega^\bullet_M$ gets identified with
 the algebra, $\f_\cM$, of smooth functions on $\cM$, so that the Lie
 algebra of derivations
 of  $\Omega^\bullet_M$ gets identified with the Lie algebra, $T_\cM$, of smooth
 vector fields on $\cM$. The de Rham differential becomes a degree one vector field $d$
 on $\cM$ satisfying $[d,d]=0$.

 \sip

There is a natural fibration $\pi: \cM \rar M$ so that we can define a subset
 $V\subset T_\cM$ by the exact sequence,
$$
0\lon V \lon T_\cM \lon \pi^*T_M \lon 0.
$$
This is the space of vertical vector fields on $\cM$ and hence it corresponds precisely
to the space of those derivations of the de Rham algebra
$\Omega^\bullet_M$ which vanish on functions,
$\Omega^0_M=\f_M$,
i.e.\ which are {\em tensorial}\,  maps $\Omega^1_M\rar  \Omega^\bullet_M$. Thus
$V$ is canonically isomorphic to the vector space of global sections of the bundle
$T_M\ot \Omega^\bullet_M$ (which we denote by the same symbol  as the bundle).
Let $i: T_M\ot \Omega^\bullet_M \rar V$ be that isomorphism. As $V$ is an integrable
distribution, it is a Lie subalgebra of $T_\cM$. Hence the standard Lie algebra structure
on $T_\cM$ induces upon restriction to $V$ and the isomorphism $i$
 a Lie algebra structure, $[\, ,\, ]_{NR}$, on $T_M\ot \Omega^\bullet_M$
which is often called in the
literature the {\em Nijenhuis-Richardson brackets} (see, e.g., \cite{Mi}).

\sip

Consider a map,
$$
\Ba{rccc}
\circ: &  T_M\ot \Omega^\bullet_M \times T_M\ot \Omega^\bullet_M   & \lon &
T_M\ot \Omega^\bullet_M \\
       &  v_1\ot \omega_1 \times   v_2\ot \omega_2 \lon &  v_2\ot (\om_1\wedge
       v_1\lrcorner \omega_2).
\Ea
$$
It is easy to check that
$$
[a,b]_{NR}=a\circ b - (-1)^{(p+1)(q+1)}b\circ a
$$
for any $a\in T_M\ot \Omega^p_M$ and  $b\in T_M\ot \Omega^q_M$.

\sip

Consider next a subspace of $T_\cM$ consisting of vector fields which commute with
the de Rham vector field $d$.
It is in fact canonically isomorphic to $V[-1]$: if $X\in T_\cM$
satisfies $[d,X]=0$ then $X=[d, i(Y)]$ for some uniquely defined
$Y\in T_M\ot \Omega^\bullet_M$. Moreover, this subspace is again a Lie subalgebra
of $T_\cM$. Hence the standard Lie algebra structure on $T_\cM$ induces on
$V\simeq T_M\ot \Omega^\bullet_M$ the structure, $[\, \bullet\, ]$, of odd Lie algebra,
$$
\left[[d,i(X_1)], [d,i(X_2)]\right]=:\left[d,i([X_1\bullet X_2])\right].
$$
These brackets $[\, \bullet\, ]$ are often called in the literature
the {\em Fr\"olicher-Nijenhuis brackets}.

\sip

An arbitrary vector field, $X$, on $\cM$ can always be decomposed into a sum,
$$
X = i(X_1) + [d, i(X_2)],
$$
for some uniquely defined $X_1, X_2\in  T_M\ot \Omega^\bullet_M$. Put another
way, $T_\cM= V\oplus V[-1]$ and, as it is not hard to check, the
standard Lie bracket on $T_\cM$ gets represented
in terms of operations $\circ$ and $[\, \bullet\, ]$ precisely as in
Section 3.1.2.

\sip

We conclude that $ T_M\ot \Omega^\bullet_M$ together with
the  Nijenhuis-Richardson composition $\circ$ and the Fr\"olicher-Nijenhuis bracket
$[\, \bullet\, ]$ is a pre-Lie$^{\mathbf\frak 2}$ algebra.

\bip

{\bf 3.2. Operads, free operads and trees.}
An operad is a pair of collections,
 $$
 \f=\left( \left\{\f(n)\right\}_{n\geq 1}\
 ,\  \{\circ^{n,n'}_{i}\}_{ n,n'\geq 1 \atop 1\leq i\leq n}\right),
 $$
 where
\Bi
\item[(i)]
 $\left\{\f(n)\right\}_{n\geq 1}$ is an an $\Bbb S$-{module}, i.e.\
 a collection of $\Z$-graded vector spaces $\f(n)$
 equipped with a linear right action of the
 permutation group $\bS_n$;
\item[(ii)] $\{\circ^{n,n'}_{i}\}_{ n,n'\geq 1 \atop 1\leq i\leq n}$ is a collection
of  linear equivariant
 maps,
 $\circ_i^{n,n'}: \f(n)\ot \f(n') \lon \f(n+n'-1)$,
 satisfying associativity conditions which can be found, e.g., in \cite{MSS}.
 \Ei

\sip
A differential graded (shortly, dg) operad  is an operad $\f$
equipped
 with a degree 1 equivariant linear map $d: \f(n) \rar \f(n)$, $\forall n$,
 satisfying the conditions,
 \Beqrn
 d^2 &=& 0, \\
 d\left(f\circ_i^{n,n'} f'\right) &=&
 (df)\circ_i^{n,n'} f' + (-1)^{|f|} f\circ_i^{n,n'} df', \ \  \forall f\in \f(n), f'\in \f(n').
 \Eeqrn
 The associated cohomology $\bS$-module $H(\f):= \{H^{\bullet}(\f(n))\}_{n\geq 1}$
 has an induced operad structure.

 \sip

 An {\em ideal}\, in an operad $\f$ is a collection $I$
 of $\bS_n$-invariant
 subspaces $\{\cI(n)\subset \f(n)\}_{n\geq 1}$
 such that $f\circ_i^{n,n'} f'\in \cI(n+n'-1)$
 whenever $f\in \cI(n)$ or $f'\in \cI(n')$; in particular, $\cI$ is a suboperad of $\f$.
 It is clear that the quotient ${\Bbb S}$-module
 $\{\f(n)/ \cI(n)\}_{n\geq 1}$ has an induced structure of operad called the
 {\em quotient operad}.

 \sip

An operad $\f$ with $\f(1)=0$ is called {\em simply connected}.
\sip

With any graded vector space $V$ one can associate an {\em endomorphism}\, operad,
 $$
 \cE_V=\left\{ \cE_V(n):= Hom(V^{\ot n}, V),
  \{\circ^{n,n'}_{i}\}_{ n,n'\geq 1 \atop 1\leq i\leq n}\right\}
 $$
where the composition,
 $f\circ_i^{n,n'} f'$, is  just insertion of the output of $f'\in Hom(V^{\ot n'},V)$
 into the $i$-th input of $f\in Hom(V^{\ot n},V)$.
\sip

 An {\em algebra over an operad}\,  $\f$ is, by definition, a $\Z$-graded vector space $V$
 together with
 a morphism of operads $\f\rar \cE_V$.

\sip

With any $\bS$-module  $\cE=\{\cE(n)\}_{n\geq 1}$ one can  associate an operad,
$Free(\cE)$, called the {\em free operad}\,
 on $\cE$. Its construction uses decorated trees
defined as follows.

\sip

Let $n$ be a natural number and $[n]$ denote the set $\{1,2,\ldots, n\}$.
 An {\em $[n]$-tree}\, $T$ is, by definition, the data $(V_T, N_T, \phi_T)$
 consisting of
 \Bi
 \item a stratified finite set $V_T= V^i_T \sqcup V^t_T$ whose elements are called
 {\em vertices};
  elements of the subset $V^i_T$ (resp.\ $V^t_T$) are called {\em internal}\, (resp.\ {\em tail})
 { vetices};
 \item a bijection $\phi: V^t_T \rar [n]$;
 \item a map $N_T: V_T \rar V_T$ satisfying the conditions: (i) $N_T$ has
 only one fixed point
 $root_T$ which lies in $V^i_T$ and is called the {\em root vertex}, (ii)
 $N_T^k(v)=root_T$,
 $\forall\ v\in V_T$ and $k\gg 1$, (iii) for all $v\in V_T^i$ the cardinality, $\# v$,
 of the set $N_T^{-1}(v)$ is greater than or equal to $1$, while for all $v\in V_T^t$ one has
  $\# v=0$.
 \Ei

 The number $\# v$ is often called the {\em valency}\, of the vertex $v$; the pairs
 $(v, N_T(v))$ are called edges.

 \sip

 Given an $\bS$-module $\cE=\{\cE(n)\}_{n\geq 1}$, we
 can associate to an $[n]$-tree $T$ the vector space
 $$
 \cE(T):= \bigotimes_{v\in V_T^i} \cE(\# v).
 $$
 Its elements are interpreted as $[n]$-trees
 whose internal vertices are decorated with elements of $\cE$.
 The permutation group $\bS_n$  acts on this space via
 relabelling the tail vertices (i.e\ changing
 $\phi_T$ to $\sigma\circ \phi_T$, $\sigma\in \bS_n$).

 \sip

 Now, as an $\bS$-module the free operad $Free(\cE)$ is defined as
 $$
 Free(\cE)(n) = \bigoplus_{[[n]-{\mathrm trees}\ T]} \cE(T),
 $$
 where the summation goes over all isomorphism classes of $[n]$-trees.
 The composition,
 say $f\circ_i^{n,n'} f'$, is given by gluing the root vertex of the
 decorated $[n]$-tree
 $f\in Free(\cE)(n)$  with the $i$-labelled
 tail vertex of the decorated $[n']$-tree $f'$. The new numeration,
 $\phi: V_T^t\rar [n+n'-1]$,
  of tails is clear.

\sip

For an $\bS$-module $\f=\{\f(n)\}_{n\geq 1}$
 we define $\f\{m\}$ to be the $\bS$-module given by the tensor product,
 $$
 \f\{m\}(n):= \f(n)\ot_k sgn_n^{\ot m}[m(1-n)],
 $$
 where $sgn_n$ is the sign representation of the permutation group $\bS_n$. If $\f$
 is a dg operad, then $\f\{m\}$ is naturally a dg operad as well:
  structure of $\f\{m\}$-algebra on a dg vector space $V$ is the same as structure
 of $\f$-algebra on the shifted dg vector space $V[m]$.

 \sip

 Let $\f=\{\f(n),\circ_i^{n,n'}, d\}$ be a simply connected operad with all
 vector spaces $\f(n)$ being finite dimensional,
 and let $\f^*[-1]$ stand for the
 $\bS$-module $\{\f(n)^*[-1]\}_{n\geq 2}$. It was shown in \cite{GK} that the free operad
 $Free(\f^*[-1]\{-1\})$ can be naturally made into a {\em differential}\,  operad,
 $$
 {\bf D}\f=(Free(\f^*[-1]\{-1\}), d),
 $$
 with the differential $d$ defined by the dualization of the
 operadic compositions $\circ_i^{n,n'}$. This  construction plays a
 special role in the theory of so-called {\em quadratic}\,  operads,
 the ones which  can be represented as a quotient,
 $$
 \f=\frac{Free(\cE)}{<R>},
 $$
 of the free operad generated by an $\bS$-module $\cE$ with $\cE(n)=0$ for $n\neq 2$
 by an ideal generated by an $\bS_3$-invariant subspace $R$ in
 $ Free(\cE)(3)$.

 \sip

 The {\em Koszul dual}\, of a quadratic operad $\f=Free(\cE)/<R>$ is, by definition,
 the quadratic operad $\f^!=Free(\check{\cE})/<R^{\bot}>$ where $\check{\cE}$ is the
 $\bS$-module
 whose only non-vanishing component is $\check{\cE}(2)=\cE(2)^*\ot sgn
 _2$
 and $R^{\bot}$ is the annihilator of $R$, i.e.\
 the kernel of the natural map $Free(\check{\cE})(3) \rar R^*$.

 \sip

 Applying the cobar construction to the Koszul dual of a quadratic operad $\f$ one gets
 a cofibrant  dg operad ${\bf D}\f^!$ together with  a canonical map of dg operads \cite{GK},
 $$
 ({\bf D}\f^!, d) \lon (\f,0).
 $$
If this map happens to be an isomorphism, the operad $\f$ is called {\em Koszul}
and then $ ({\bf D}\f^!, d)$ provides the minimal dg resolution of $\f$.

\sip

For any simply connected operad $\f$  the associated  operad ${\bf D(\f^!)}$ is minimal,
i.e.\ its differential is decomposable.
Thus the
 {\em cobar construction}, that is the functor,
$$
\f \rar \f_\infty:= ({\bf D}\f^!, d),
$$
associates with {\em any}\,
 quadratic operad a minimal dg operad whose algebras are strongly homotopy
 ones (i.e.\ are transferable via quasi-isomorphisms of complexes) and whose
 homology operad
 covers $\f$.

 \bip

{\bf 3.3. Operad of pre-Lie$^{\bf\frak 2}$ algebras.} Let us denote this operad by
$\cP$. Clearly, it is quadratic, $\cP={Free(\cE)}/{<R>}$, with
$\cE(2):=k[\bS_2][0]\oplus {\mathbf 1}_2[-1]$ and the relations 3.1(i)-(iii).
Here and elsewhere
${\mathbf 1}_n$ (respectively, $k[\bS_n]$)
stands for the one dimensional trivial (respectively, $n!$ dimensional regular)
  representation of $\bS_n$.

\sip

As the free operad $Free(\cE)$
 is a construction on binary trees, it will be useful
 to have corolla type notation for the generators of $\cE(2)$: we shall
 represent the two basis vectors, $(1)(2)$ and $(12)$, of the summand
$k[\bS_2][0]$ by the labelled planar corollas of degree 0,
 $$
 
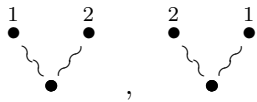
\begin{xy}
 <5mm,0cm>*{\bullet};<0cm,7mm>*{\bullet}**@{~},
 <5mm,0cm>*{\bullet};<10mm,7mm>*{\bullet}**@{~},
 <5mm,0cm>*{\bullet};<0mm,9mm>*{^1}**@{},
 <5mm,0cm>*{\bullet};<10mm,9mm>*{^2}**@{},
 \end{xy}
 \ \ \ , \ \ \
 \begin{xy}
 <5mm,0cm>*{\bullet};<0cm,7mm>*{\bullet}**@{~},
 <5mm,0cm>*{\bullet};<10mm,7mm>*{\bullet}**@{~},
 <5mm,0cm>*{\bullet};<0mm,9mm>*{^2}**@{},
 <5mm,0cm>*{\bullet};<10mm,9mm>*{^1}**@{},
 \end{xy} \ \  \ \
 $$
while the basis vector of the summand ${\mathbf 1}_2[-1]$  by the
unique space corolla of degree 1,
$$
\begin{xy}
 <5mm,0cm>*{\bullet};<0cm,7mm>*{\bullet}**@{-},
 <5mm,0cm>*{\bullet};<10mm,7mm>*{\bullet}**@{-},
 <5mm,0cm>*{\bullet};<0mm,9mm>*{^1}**@{},
 <5mm,0cm>*{\bullet};<10mm,9mm>*{^2}**@{},
 \end{xy}
 \ \ =
\ \
\begin{xy}
 <5mm,0cm>*{\bullet};<0cm,7mm>*{\bullet}**@{-},
 <5mm,0cm>*{\bullet};<10mm,7mm>*{\bullet}**@{-},
 <5mm,0cm>*{\bullet};<0mm,9mm>*{^2}**@{},
 <5mm,0cm>*{\bullet};<10mm,9mm>*{^1}**@{},
 \end{xy}
$$
Then the ideal ${<R>}$ is generated by the following elements of of $Free(\cE)(3)$,
$$
 \begin{xy}
 <5mm,0mm>*{\bullet};<0cm,7mm>*{\bullet}**@{~},
 <5mm,0cm>*{\bullet};<10mm,7mm>*{\bullet}**@{~},
 <0mm,7mm>*{\bullet};<-5mm,14mm>*{\bullet}**@{~},
 <0mm,7mm>*{\bullet};<5mm,14mm>*{\bullet}**@{~},
 <5mm,0cm>*{\bullet};<-5mm,16.5mm>*{^{i_1}}**@{},
 <5mm,0cm>*{\bullet};<5mm,16.5mm>*{^{i_2}}**@{},
 <5mm,0cm>*{\bullet};<11mm,9.3mm>*{^{i_3}}**@{},
 \end{xy}
 \
 -
 \
  \begin{xy}
 <5mm,0mm>*{\bullet};<0cm,7mm>*{\bullet}**@{~},
 <5mm,0cm>*{\bullet};<10mm,7mm>*{\bullet}**@{~},
 <0mm,7mm>*{\bullet};<-5mm,14mm>*{\bullet}**@{~},
 <0mm,7mm>*{\bullet};<5mm,14mm>*{\bullet}**@{~},
 <5mm,0cm>*{\bullet};<-5mm,16.5mm>*{^{i_1}}**@{},
 <5mm,0cm>*{\bullet};<5mm,16.5mm>*{^{i_3}}**@{},
 <5mm,0cm>*{\bullet};<11mm,9.3mm>*{^{i_2}}**@{},
 \end{xy}
 \
 -
 \
 \begin{xy}
 <5mm,0mm>*{\bullet};<0cm,7mm>*{\bullet}**@{~},
 <5mm,0cm>*{\bullet};<10mm,7mm>*{\bullet}**@{~},
 <10mm,7mm>*{\bullet};<5mm,14mm>*{\bullet}**@{~},
 <10mm,7mm>*{\bullet};<15mm,14mm>*{\bullet}**@{~},
 <5mm,0cm>*{\bullet};<-1mm,9.3mm>*{^{i_1}}**@{},
 <5mm,0cm>*{\bullet};<5mm,16.5mm>*{^{i_2}}**@{},
 <5mm,0cm>*{\bullet};<16mm,16.5mm>*{^{i_3}}**@{},
 \end{xy}
\
+
\
\begin{xy}
 <5mm,0mm>*{\bullet};<0cm,7mm>*{\bullet}**@{~},
 <5mm,0cm>*{\bullet};<10mm,7mm>*{\bullet}**@{~},
 <10mm,7mm>*{\bullet};<5mm,14mm>*{\bullet}**@{~},
 <10mm,7mm>*{\bullet};<15mm,14mm>*{\bullet}**@{~},
 <5mm,0cm>*{\bullet};<-1mm,9.3mm>*{^{i_1}}**@{},
 <5mm,0cm>*{\bullet};<5mm,16.5mm>*{^{i_3}}**@{},
 <5mm,0cm>*{\bullet};<16mm,16.5mm>*{^{i_2}}**@{},
 \end{xy},
 $$
\begin{equation}\label{A}
 \begin{xy}
 <5mm,0mm>*{\bullet};<0cm,7mm>*{\bullet}**@{-},
 <5mm,0cm>*{\bullet};<10mm,7mm>*{\bullet}**@{-},
 <0mm,7mm>*{\bullet};<-5mm,14mm>*{\bullet}**@{-},
 <0mm,7mm>*{\bullet};<5mm,14mm>*{\bullet}**@{-},
 <5mm,0cm>*{\bullet};<-5.4mm,16.3mm>*{^1}**@{},
 <5mm,0cm>*{\bullet};<5mm,16.3mm>*{^2}**@{},
 <5mm,0cm>*{\bullet};<10.6mm,9.1mm>*{^3}**@{},
 \end{xy}\ \ + \ \
\begin{xy}
 <5mm,0mm>*{\bullet};<0cm,7mm>*{\bullet}**@{-},
 <5mm,0cm>*{\bullet};<10mm,7mm>*{\bullet}**@{-},
 <0mm,7mm>*{\bullet};<-5mm,14mm>*{\bullet}**@{-},
 <0mm,7mm>*{\bullet};<5mm,14mm>*{\bullet}**@{-},
 <5mm,0cm>*{\bullet};<-5.4mm,16.3mm>*{^3}**@{},
 <5mm,0cm>*{\bullet};<5mm,16.3mm>*{^1}**@{},
 <5mm,0cm>*{\bullet};<10.6mm,9.1mm>*{^2}**@{},
 \end{xy}\ \ + \ \
\begin{xy}
 <5mm,0mm>*{\bullet};<0cm,7mm>*{\bullet}**@{-},
 <5mm,0cm>*{\bullet};<10mm,7mm>*{\bullet}**@{-},
 <0mm,7mm>*{\bullet};<-5mm,14mm>*{\bullet}**@{-},
 <0mm,7mm>*{\bullet};<5mm,14mm>*{\bullet}**@{-},
 <5mm,0cm>*{\bullet};<-5.4mm,16.3mm>*{^2}**@{},
 <5mm,0cm>*{\bullet};<5mm,16.3mm>*{^3}**@{},
 <5mm,0cm>*{\bullet};<10.6mm,9.1mm>*{^1}**@{},
 \end{xy}\ ,
\end{equation}
 \begin{equation}\label{B}
 \begin{xy}
 <5mm,0mm>*{\bullet};<0cm,7mm>*{\bullet}**@{~},
 <5mm,0cm>*{\bullet};<10mm,7mm>*{\bullet}**@{~},
 <0mm,7mm>*{\bullet};<-5mm,14mm>*{\bullet}**@{-},
 <0mm,7mm>*{\bullet};<5mm,14mm>*{\bullet}**@{-},
 <5mm,0cm>*{\bullet};<-5mm,16.4mm>*{^{i_1}}**@{},
 <5mm,0cm>*{\bullet};<5mm,16.4mm>*{^{i_2}}**@{},
 <5mm,0cm>*{\bullet};<11mm,9.3mm>*{^{i_3}}**@{},
  \end{xy}\  + \
  \begin{xy}
 <5mm,0mm>*{\bullet};<0cm,7mm>*{\bullet}**@{~},
 <5mm,0cm>*{\bullet};<10mm,7mm>*{\bullet}**@{~},
 <10mm,7mm>*{\bullet};<5mm,14mm>*{\bullet}**@{-},
 <10mm,7mm>*{\bullet};<15mm,14mm>*{\bullet}**@{-},
 <5mm,0cm>*{\bullet};<-1mm,9.3mm>*{^{i_1}}**@{},
 <5mm,0cm>*{\bullet};<5mm,16.4mm>*{^{i_2}}**@{},
 <5mm,0cm>*{\bullet};<16mm,16.4mm>*{^{i_3}}**@{},
 \end{xy}
  \ + \
 \begin{xy}
 <5mm,0mm>*{\bullet};<0cm,7mm>*{\bullet}**@{~},
 <5mm,0cm>*{\bullet};<10mm,7mm>*{\bullet}**@{~},
 <10mm,7mm>*{\bullet};<5mm,14mm>*{\bullet}**@{-},
 <10mm,7mm>*{\bullet};<15mm,14mm>*{\bullet}**@{-},
 <5mm,0cm>*{\bullet};<-1mm,9.3mm>*{^{i_2}}**@{},
 <5mm,0cm>*{\bullet};<5mm,16.4mm>*{^{i_1}}**@{},
 <5mm,0cm>*{\bullet};<16mm,16.4mm>*{^{i_3}}**@{},
 \end{xy}
 \ - \
 \begin{xy}
 <5mm,0mm>*{\bullet};<0cm,7mm>*{\bullet}**@{-},
 <5mm,0cm>*{\bullet};<10mm,7mm>*{\bullet}**@{-},
 <10mm,7mm>*{\bullet};<5mm,14mm>*{\bullet}**@{~},
 <10mm,7mm>*{\bullet};<15mm,14mm>*{\bullet}**@{~},
 <5mm,0cm>*{\bullet};<-1mm,9.3mm>*{^{i_1}}**@{},
 <5mm,0cm>*{\bullet};<5mm,16.4mm>*{^{i_2}}**@{},
 <5mm,0cm>*{\bullet};<16mm,16.4mm>*{^{i_3}}**@{},
 \end{xy}
   \ -\
\begin{xy}
 <5mm,0mm>*{\bullet};<0cm,7mm>*{\bullet}**@{-},
 <5mm,0cm>*{\bullet};<10mm,7mm>*{\bullet}**@{-},
 <10mm,7mm>*{\bullet};<5mm,14mm>*{\bullet}**@{~},
 <10mm,7mm>*{\bullet};<15mm,14mm>*{\bullet}**@{~},
 <5mm,0cm>*{\bullet};<-1mm,9.3mm>*{^{i_2}}**@{},
 <5mm,0cm>*{\bullet};<5mm,16.4mm>*{^{i_1}}**@{},
 <5mm,0cm>*{\bullet};<16mm,16.4mm>*{^{i_3}}**@{},
 \end{xy},
 \end{equation}
where $(i_1,i_2,i_3)$ is an arbitrary permutation of $(1,2,3)$.

 \bip

{\bf 3.4. Koszul dual of $\cP$.} By definition,
$$
\check{\cE}(2)=\cE(2)^*\ot sgn _2= k[\bS_2][0]\oplus sgn_2[1].
 $$
If we represent the basis vectors of $\check{\cE}(2)$ by two planar
 corollas of degree 0,
$$
 \begin{xy}
 <5mm,0cm>*{\bullet};<0cm,7mm>*{\bullet}**@{~},
 <5mm,0cm>*{\bullet};<10mm,7mm>*{\bullet}**@{~},
 <5mm,0cm>*{\bullet};<-0.2mm,9mm>*{^1}**@{},
 <5mm,0cm>*{\bullet};<10mm,9mm>*{^2}**@{},
 \end{xy}
 \ \ \ , \ \ \
 \begin{xy}
 <5mm,0cm>*{\bullet};<0cm,7mm>*{\bullet}**@{~},
 <5mm,0cm>*{\bullet};<10mm,7mm>*{\bullet}**@{~},
 <5mm,0cm>*{\bullet};<-0mm,9mm>*{^2}**@{},
 <5mm,0cm>*{\bullet};<10mm,9mm>*{^1}**@{},
 \end{xy} \ \  \ \
 $$
and one corolla of degree $-1$,
$$
\begin{xy}
 <5mm,0cm>*{\bullet};<0cm,7mm>*{\bullet}**@{-},
 <5mm,0cm>*{\bullet};<10mm,7mm>*{\bullet}**@{-},
 <5mm,0cm>*{\bullet};<-0.2mm,9mm>*{^1}**@{},
 <5mm,0cm>*{\bullet};<10mm,9mm>*{^2}**@{},
 \end{xy}
 \ \ =
\ - \
\begin{xy}
 <5mm,0cm>*{\bullet};<0cm,7mm>*{\bullet}**@{-},
 <5mm,0cm>*{\bullet};<10mm,7mm>*{\bullet}**@{-},
 <5mm,0cm>*{\bullet};<0mm,9mm>*{^2}**@{},
 <5mm,0cm>*{\bullet};<10mm,9mm>*{^1}**@{},
 \end{xy}
$$
then the Koszul dual operad $\cP^!$ is spanned by all possible binary trees on these
corollas, modulo the following relations, i.e.\ generators
of  $R^{\bot}\subset Free(\hat{\cE})(3)$,
\begin{equation}\label{C}
 \begin{xy}
 <5mm,0mm>*{\bullet};<0cm,7mm>*{\bullet}**@{~},
 <5mm,0cm>*{\bullet};<10mm,7mm>*{\bullet}**@{~},
 <0mm,7mm>*{\bullet};<-5mm,14mm>*{\bullet}**@{~},
 <0mm,7mm>*{\bullet};<5mm,14mm>*{\bullet}**@{~},
 <5mm,0cm>*{\bullet};<-5mm,16.4mm>*{^{i_1}}**@{},
 <5mm,0cm>*{\bullet};<5mm,16.4mm>*{^{i_2}}**@{},
 <5mm,0cm>*{\bullet};<11mm,9.2mm>*{^{i_3}}**@{},
 \end{xy}
 \
 -
 \
 \begin{xy}
 <5mm,0mm>*{\bullet};<0cm,7mm>*{\bullet}**@{~},
 <5mm,0cm>*{\bullet};<10mm,7mm>*{\bullet}**@{~},
 <10mm,7mm>*{\bullet};<5mm,14mm>*{\bullet}**@{~},
 <10mm,7mm>*{\bullet};<15mm,14mm>*{\bullet}**@{~},
 <5mm,0cm>*{\bullet};<-1mm,9.2mm>*{^{i_1}}**@{},
 <5mm,0cm>*{\bullet};<5mm,16.4mm>*{^{i_2}}**@{},
 <5mm,0cm>*{\bullet};<16mm,16.4mm>*{^{i_3}}**@{},
 \end{xy}
 \ = \
 \begin{xy}
 <5mm,0mm>*{\bullet};<0cm,7mm>*{\bullet}**@{-},
 <5mm,0cm>*{\bullet};<10mm,7mm>*{\bullet}**@{-},
 <0mm,7mm>*{\bullet};<-5mm,14mm>*{\bullet}**@{-},
 <0mm,7mm>*{\bullet};<5mm,14mm>*{\bullet}**@{-},
 <5mm,0cm>*{\bullet};<-5mm,16.4mm>*{^{i_1}}**@{},
 <5mm,0cm>*{\bullet};<5mm,16.4mm>*{^{i_2}}**@{},
 <5mm,0cm>*{\bullet};<11mm,9.2mm>*{^{i_3}}**@{},
 \end{xy}
 \
 -
 \
 \begin{xy}
 <5mm,0mm>*{\bullet};<0cm,7mm>*{\bullet}**@{-},
 <5mm,0cm>*{\bullet};<10mm,7mm>*{\bullet}**@{-},
 <10mm,7mm>*{\bullet};<5mm,14mm>*{\bullet}**@{-},
 <10mm,7mm>*{\bullet};<15mm,14mm>*{\bullet}**@{-},
 <5mm,0cm>*{\bullet};<-1mm,9.2mm>*{^{i_1}}**@{},
 <5mm,0cm>*{\bullet};<5mm,16.4mm>*{^{i_2}}**@{},
 <5mm,0cm>*{\bullet};<16mm,16.4mm>*{^{i_3}}**@{},
 \end{xy}
 \ =0,
\end{equation}
\begin{equation}\label{D}
 \begin{xy}
 <5mm,0mm>*{\bullet};<0cm,7mm>*{\bullet}**@{~},
 <5mm,0cm>*{\bullet};<10mm,7mm>*{\bullet}**@{~},
 <0mm,7mm>*{\bullet};<-5mm,14mm>*{\bullet}**@{~},
 <0mm,7mm>*{\bullet};<5mm,14mm>*{\bullet}**@{~},
 <5mm,0cm>*{\bullet};<-5mm,16.4mm>*{^{i_1}}**@{},
 <5mm,0cm>*{\bullet};<5mm,16.4mm>*{^{i_2}}**@{},
 <5mm,0cm>*{\bullet};<11mm,9.2mm>*{^{i_3}}**@{},
 \end{xy}
 \ =\
\begin{xy}
 <5mm,0mm>*{\bullet};<0cm,7mm>*{\bullet}**@{~},
 <5mm,0cm>*{\bullet};<10mm,7mm>*{\bullet}**@{~},
 <0mm,7mm>*{\bullet};<-5mm,14mm>*{\bullet}**@{~},
 <0mm,7mm>*{\bullet};<5mm,14mm>*{\bullet}**@{~},
 <5mm,0cm>*{\bullet};<-5mm,16.4mm>*{^{i_1}}**@{},
 <5mm,0cm>*{\bullet};<5mm,16.4mm>*{^{i_3}}**@{},
 <5mm,0cm>*{\bullet};<11mm,9.2mm>*{^{i_2}}**@{},
 \end{xy}
\ \ \ ,\ \ \
\begin{xy}
 <5mm,0mm>*{\bullet};<0cm,7mm>*{\bullet}**@{-},
 <5mm,0cm>*{\bullet};<10mm,7mm>*{\bullet}**@{-},
 <10mm,7mm>*{\bullet};<5mm,14mm>*{\bullet}**@{~},
 <10mm,7mm>*{\bullet};<15mm,14mm>*{\bullet}**@{~},
 <5mm,0cm>*{\bullet};<-1mm,9.2mm>*{^{i_1}}**@{},
 <5mm,0cm>*{\bullet};<5mm,16.4mm>*{^{i_2}}**@{},
 <5mm,0cm>*{\bullet};<16mm,16.4mm>*{^{i_3}}**@{},
 \end{xy}
\ = \
 \begin{xy}
 <5mm,0mm>*{\bullet};<0cm,7mm>*{\bullet}**@{~},
 <5mm,0cm>*{\bullet};<10mm,7mm>*{\bullet}**@{~},
 <0mm,7mm>*{\bullet};<-5mm,14mm>*{\bullet}**@{-},
 <0mm,7mm>*{\bullet};<5mm,14mm>*{\bullet}**@{-},
 <5mm,0cm>*{\bullet};<-5mm,16.4mm>*{^{i_1}}**@{},
 <5mm,0cm>*{\bullet};<5mm,16.4mm>*{^{i_2}}**@{},
 <5mm,0cm>*{\bullet};<11mm,9.2mm>*{^{i_3}}**@{},
 \end{xy}
\end{equation}
\begin{equation}\label{E}
 \begin{xy}
 <5mm,0mm>*{\bullet};<0cm,7mm>*{\bullet}**@{~},
 <5mm,0cm>*{\bullet};<10mm,7mm>*{\bullet}**@{~},
 <10mm,7mm>*{\bullet};<5mm,14mm>*{\bullet}**@{-},
 <10mm,7mm>*{\bullet};<15mm,14mm>*{\bullet}**@{-},
 <5mm,0cm>*{\bullet};<-1mm,9.2mm>*{^{i_1}}**@{},
 <5mm,0cm>*{\bullet};<5mm,16.4mm>*{^{i_2}}**@{},
 <5mm,0cm>*{\bullet};<16mm,16.4mm>*{^{i_3}}**@{},
 \end{xy}
 \ = \
 \begin{xy}
 <5mm,0mm>*{\bullet};<0cm,7mm>*{\bullet}**@{~},
 <5mm,0cm>*{\bullet};<10mm,7mm>*{\bullet}**@{~},
 <0mm,7mm>*{\bullet};<-5mm,14mm>*{\bullet}**@{-},
 <0mm,7mm>*{\bullet};<5mm,14mm>*{\bullet}**@{-},
 <5mm,0cm>*{\bullet};<-5mm,16.4mm>*{^{i_1}}**@{},
 <5mm,0cm>*{\bullet};<5mm,16.4mm>*{^{i_3}}**@{},
 <5mm,0cm>*{\bullet};<11mm,9.2mm>*{^{i_2}}**@{},
  \end{xy}\ \ - \ \
   \begin{xy}
 <5mm,0mm>*{\bullet};<0cm,7mm>*{\bullet}**@{~},
 <5mm,0cm>*{\bullet};<10mm,7mm>*{\bullet}**@{~},
 <0mm,7mm>*{\bullet};<-5mm,14mm>*{\bullet}**@{-},
 <0mm,7mm>*{\bullet};<5mm,14mm>*{\bullet}**@{-},
 <5mm,0cm>*{\bullet};<-5mm,16.4mm>*{^{i_1}}**@{},
 <5mm,0cm>*{\bullet};<5mm,16.4mm>*{^{i_2}}**@{},
 <5mm,0cm>*{\bullet};<11mm,9.2mm>*{^{i_3}}**@{},
  \end{xy}
 \end{equation}

Thus a $\cP^!$-algebra is  a graded vector space $V$ equipped with two
associative products,
$$
\Ba{rccc}
\circ: &  V\ot V & \lon & V \\
       & a\ot b    & \lon & a\circ b
\Ea
\ \ \ \ \  , \ \ \ \
\Ba{rccc}
\bullet: & \wedge^2 V& \lon & V[-1] \\
       & a\ot b    & \lon & a\bullet b
\Ea
$$
such that
$$
(a\circ b)\circ c = (-1)^{|b||c|} (a\circ c)\circ b, \ \ \ a\bullet (b\circ c) =
(a\bullet b)\circ c,
$$
$$
a\circ (b\bullet c)= (-1)^{|b|+1} (a\bullet b)\circ c +  (-1)^{|b|(|c|+1)}
(a\bullet c)\circ b,
$$
for any $a,b,c\in V$.

\bip

{\bf 3.4.1. Example.} Let $(A, \cdot, d)$ be an arbitrary
differential graded commutative associative algebra (say, the de Rham algebra
of a smooth manifold, cf.\ Example 3.1.4). Then $A[-1]$
has a natural structure of $\cP^!$-algebra. If we denote elements of $A[-1]$
by $\Pi\al$, where $\al$ is an element of $A$ and $\Pi$ a formal symbol of degree 1,
then following two  products in $A[-1]$,
$$
\Pi \al \circ \Pi\be = \Pi (\al\cdot d\be),
$$
$$
\Pi \al \bullet \Pi\be = \Pi (\al\cdot \be),
$$
are associative and satisfy all the necessary conditions of 3.4.

\bip

{\bf 3.5. Homology of pre-Lie$^{\frak 2}$ algebras.} Let $\f$ be a quadratic operad
and $V$ an $\f$-algebra. The $\f$-algebra structure on $V$ gives rise
to a degree $-1$ differential $d_\f$
 on
 the following positively graded
vector space \cite{GK,MSS},
$$
C_{\f}V=\oplus_{n\geq 1} C^n_{\f}(V)=\oplus_{n\geq 1}(\f^!(n))^*\ot_{\bS_n} (V[1])^{\ot n}.
$$
The associated homology,
$H_\bullet(V):=H_\bullet(C_{\f}V,d_\f)$, is called the $\f$-{\em algebra homology}\,
 of $V$.

\sip

For example, if $\f$ is the operad of Lie algebras and $\fg$ a Lie algebra,
then the operadic complex $(C_{Lie}(\fg), d_{Lie})$ is precisely the
 Chevalley-Eilenberg complex,
$$
\left(\odot^\bullet(\fg[1])=\oplus_n (\wedge^{n}\fg)[n], d_{CE}\right)
$$
 which computes
the homology of the Lie algebra $\fg$ with coefficients in the trivial module.

\sip

We shall show now that the same Chevalley-Eilenberg  construction can be used to
compute the pre-Lie$^{\frak 2}$-algebra homology.

\sip

Let $\cP$ be the operad of pre-Lie$^{\frak 2}$ algebras,
and $V$ a $\cP$-algebra. According to 3.1.2, the vector space $V\oplus V[-1]$
is naturally a Lie algebra so that the vector space
$$
\odot^\bullet\left((V\oplus V[-1])[1]\right) =
\odot^\bullet\left(V[1]\oplus V\right)
$$
has the Chevalley-Eilenberg  differential $d_{CE}$.

\sip

The subspace $V\subset V\oplus V[-1]$ is a Lie subalgebra with the bracket
$[a,b]= a\circ b - (-1)^{|a||b|}b\circ a$. Hence we get one more
Chevalley-Eilenberg complex, $(\odot^\bullet (V[1]), d_{CE})$, which is a subcomplex
of the one above. Define the quotient complex by the exact sequence,
$$
0\lon (\odot^\bullet (V[1]), d_{CE}) \lon \left(\odot^\bullet(V[1]\oplus V), d_{CE}
\right) \lon (J^\bullet, d).
$$

\bip

{\bf 3.5.1. Proposition.} {\em The quotient complex $(J^\bullet, d)$
is isomorphic to the operadic complex $(C_{\cP}^\bullet V,d_\cP)$}.

\sip

Proof is a straightforward calculation.

\sip

For example, if $V$ is the  pre-Lie$^{\frak 2}$ algebra associated with the free Lie algebra
$Free(W\oplus W[-1])$  as in Example 3.1.3, then $H_1(V)=W$
and $H_2(V)=\Img Q/[\Img Q,\Img Q]\neq 0$.

\bip

\sip

{\bf 3.6. Cobar construction.} Let $\cP$ be the operad of pre-Lie$^{\frak 2}$ algebras.
The associated cobar construction,
$\cP_\infty = {\mathbf D}\cP^!$, is the free operad on the $\bS$-module,
$$
\cP^!(n)^*\ot sgn_n [n-2]=\bigoplus_{p=0}^{n-1}
{\mathrm Ind}^{\bS_n}_{\bS_{n-p}\times \bS_{p}}
{\mathbf 1}_{n-p}\ot sgn_{p}[p-1].
$$
Let us identify the basis of the $\binom{n}{p}$-dimensional summand
${\mathrm Ind}^{\bS_n}_{\bS_{n-p}\times \bS_{p}}
{\mathbf 1}_{n-p}\ot sgn_{p}[p-1]$  with planar
$[n]$-corollas of the form,
 $$
 
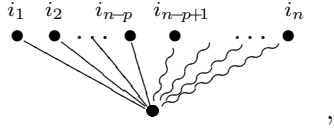
\begin{xy}
 <18mm,0cm>*{{\bullet}};<0cm,1cm>*{\bullet}**@{-},
 <18mm,0cm>*{\bullet};<5mm,1cm>*{\bullet}**@{-},
 <18mm,0cm>*{\bullet};<10mm,1cm>*{\ldots}**@{-},
 <18mm,0cm>*{\bullet};<15mm,1cm>*{\bullet}**@{-},
 <18mm,0cm>*{\bullet};<21mm,1cm>*{\bullet}**@{~},
 <18mm,0cm>*{\bullet};<26mm,1cm>*{}**@{~},
 <18mm,0cm>*{\bullet};<31mm,1cm>*{\ldots}**@{~},
 <18mm,0cm>*{\bullet};<36mm,1cm>*{\bullet}**@{~},
 <18mm,0cm>*{\bullet};<0mm,13mm>*{^{i_1}}**@{},
 <18mm,0cm>*{\bullet};<5mm,13mm>*{^{i_2}}**@{},
 <18mm,0cm>*{\bullet};<13mm,13mm>*{^{i_{n\hspace{-0.5mm}-\hspace{-0.5mm}p}}}**@{},
 <18mm,0cm>*{\bullet};<22mm,13mm>*{^{i_{n\hspace{-0.5mm}-\hspace{-0.5mm}p\hspace{-0.4mm}
 +\hspace{-0.5mm}1}}}**@{},
 <18mm,0cm>*{\bullet};<37mm,13mm>*{^{i_n}}**@{},
 \end{xy} \ \ ,
 $$
which are symmetric over the first $n-p$ inputs and antisymmetric over
the last $p$ ones, and has degree $1-p$.
As an $\bS$-module the operad $\cP_\infty$
 is the linear span
 of all possible (isomorphism classes of)  plane
 trees formed by these corollas. The compositions $\circ_i^{n,n'}$ in $\cP_\infty$
 are
 given  simply by gluing  the root vertex of an $[n']$-tree to the $i$th tail vertex
 of an $[n]$-tree. To complete the picture we need just the cobar differential.

 \sip

 {\bf 3.6.1. Proposition}. {\em The cobar differential in $\cP_{\infty}$
 is given on generators by}

\Beqrn
 \begin{xy}
<18mm,0cm>*{{\bullet}};<0cm,1cm>*{\bullet}**@{-},
 <18mm,0cm>*{\bullet};<5mm,1cm>*{\bullet}**@{-},
 <18mm,0cm>*{\bullet};<10mm,1cm>*{\ldots}**@{-},
 <18mm,0cm>*{\bullet};<15mm,1cm>*{\bullet}**@{-},
 <18mm,0cm>*{\bullet};<21mm,1cm>*{\bullet}**@{~},
 <18mm,0cm>*{\bullet};<26mm,1cm>*{}**@{~},
 <18mm,0cm>*{\bullet};<31mm,1cm>*{\ldots}**@{~},
 <18mm,0cm>*{\bullet};<36mm,1cm>*{\bullet}**@{~},
 <18mm,0cm>*{\bullet};<-0.5mm,12mm>*{^1}**@{},
 <18mm,0cm>*{\bullet};<4.5mm,12mm>*{^2}**@{},
 <18mm,0cm>*{\bullet};<13.5mm,12mm>*{^{{n\hspace{-0.5mm}-\hspace{-0.5mm}p}}}**@{},
 <18mm,0cm>*{\bullet};<22.5mm,12mm>*{^{{n\hspace{-0.5mm}-\hspace{-0.5mm}p\hspace{-0.4mm}
 +\hspace{-0.5mm}1}}}**@{},
 <18mm,0cm>*{{\bullet}};<-4mm,5mm>*{d}**@{},
 <18mm,0cm>*{\bullet};<36.5mm,12mm>*{^{n}}**@{},
 \end{xy}
 &=& \vspace{30mm}
 \sum_{I_1\sqcup I_2=(1,\ldots,n-p) \atop
{J_1\sqcup J_2=(n-p+1,\ldots,n) \atop
 {\#I_2\geq 1, \#I_1+\#J_2\geq 1\atop
 { \#I_2+\#J_1\geq 2}}}}
 \hspace{-5mm} (-1)^{\#J_2 + \sigma(J_1\sqcup J_2)}
\begin{xy}
 <28mm,0cm>*{{\bullet}};<7mm,1cm>*{\bullet}**@{-},
 <28mm,0cm>*{\bullet};<14mm,1cm>*{\ldots}**@{-},
 <28mm,0cm>*{\bullet};<20mm,1cm>*{\bullet}**@{-},
 <28mm,0cm>*{\bullet};<27mm,1cm>*{\bullet}**@{-},
 <28mm,0cm>*{\bullet};<34mm,1cm>*{\bullet}**@{~},
 <28mm,0cm>*{\bullet};<40mm,1cm>*{\ldots}**@{~},
 <28mm,0cm>*{\bullet};<48mm,1cm>*{\bullet}**@{~},
 <28mm,0cm>*{\bullet};<42mm,14.5mm>*{^{J_2}}**@{},
 <28mm,0cm>*{\bullet};<41mm,12mm>*{\overbrace{\ \ \ \ \ \  \ \ \ \ \ \ \ }}**@{},
 <27mm,1cm>*{{\bullet}};<14mm,22mm>*{\bullet}**@{-},
 <27mm,1cm>*{\bullet};<19mm,22mm>*{\cdots}**@{-},
 <27mm,1cm>*{\bullet};<24mm,22mm>*{\bullet}**@{-},
 <27mm,1cm>*{\bullet};<30mm,22mm>*{\bullet}**@{~},
 <27mm,1cm>*{\bullet};<35mm,22mm>*{\cdots}**@{~},
 <27mm,1cm>*{\bullet};<40mm,22mm>*{\bullet}**@{~},
 <38mm,0cm>*{};<35mm,27mm>*{^{J_1}}**@{},
 <38mm,0cm>*{};<35mm,24mm>*{\overbrace{\ \ \ \ \ \ \ \ \ \ \ }}**@{},
<38mm,0cm>*{};<19mm,27mm>*{^{I_2}}**@{},
 <38mm,0cm>*{};<19mm,24mm>*{\overbrace{\ \ \  \ \ \ \ \ \ \ \  }}**@{},
 <38mm,0cm>*{};<13mm,14.5mm>*{^{I_1}}**@{},
 <38mm,0cm>*{};<13mm,12mm>*{\overbrace{\ \ \ \ \ \ \ \ \ \ \ \ \ }}**@{},
 \end{xy}
 \\
&&
\vspace{30mm}
 -\sum_{I_1\sqcup I_2=(1,\ldots,n-p) \atop
{J_1\sqcup J_2\sqcup J_3=(n-p+1,\ldots,n) \atop
 {\#I_1\geq 1, \#I_2\geq 1\atop
 { \#I_1+\#J_3\geq 1 ,\#I_2+\#J_2\geq 1}}}}
 \hspace{-5mm} (-1)^r
  \hspace{-5mm}
\begin{xy}
 <28mm,0cm>*{{\bullet}};<7mm,1cm>*{\bullet}**@{-},
 <28mm,0cm>*{\bullet};<14mm,1cm>*{\ldots}**@{-},
 <28mm,0cm>*{\bullet};<20mm,1cm>*{\bullet}**@{-},
 <28mm,0cm>*{\bullet};<27mm,1cm>*{\bullet}**@{~},
 <28mm,0cm>*{\bullet};<34mm,1cm>*{\bullet}**@{~},
 <28mm,0cm>*{\bullet};<40mm,1cm>*{\ldots}**@{~},
 <28mm,0cm>*{\bullet};<48mm,1cm>*{\bullet}**@{~},
 <28mm,0cm>*{\bullet};<42mm,14.5mm>*{^{J_3}}**@{},
 <28mm,0cm>*{\bullet};<41mm,12mm>*{\overbrace{\ \ \ \ \ \  \ \ \ \ \ \ \ }}**@{},
 <27mm,1cm>*{{\bullet}};<11mm,22mm>*{\bullet}**@{-},
 <27mm,1cm>*{\bullet};<16mm,22mm>*{\cdots}**@{-},
 <27mm,1cm>*{\bullet};<21mm,22mm>*{\bullet}**@{-},
 <27mm,1cm>*{\bullet};<33mm,22mm>*{\bullet}**@{~},
    <27mm,1cm>*{\bullet};<27mm,22mm>*{\bullet}**@{-},
    <38mm,0cm>*{};<27mm,25mm>*{^{J_1}}**@{},
 <27mm,1cm>*{\bullet};<38mm,22mm>*{\cdots}**@{~},
 <27mm,1cm>*{\bullet};<43mm,22mm>*{\bullet}**@{~},
 <38mm,0cm>*{};<38mm,27mm>*{^{J_2}}**@{},
 <38mm,0cm>*{};<38mm,24mm>*{\overbrace{\ \ \ \ \ \ \ \ \ \ \ }}**@{},
<38mm,0cm>*{};<16mm,27mm>*{^{I_2}}**@{},
 <38mm,0cm>*{};<16mm,24mm>*{\overbrace{\ \ \  \ \ \ \ \ \ \   }}**@{},
 <38mm,0cm>*{};<13mm,14.5mm>*{^{I_1}}**@{},
 <38mm,0cm>*{};<13mm,12mm>*{\overbrace{\ \ \ \ \ \ \ \ \ \ \ \ \ }}**@{},
 \end{xy}
 \Eeqrn
{\em where $r={\#J_2 + \#J_3 +\sigma(J_1\sqcup J_2\sqcup J_3)}$ and
$\sigma(J_1\sqcup J_2)$ and ${\sigma(J_1\sqcup J_2\sqcup J_3)}$ stand for the parities
of the permutations  $(n-p+1,\ldots,n)\rar (J_1\sqcup J_2)$
and, respectively, $(n-p+1,\ldots,n)\rar (J_1\sqcup J_2\sqcup J_3)$.}

 \bip

 Proof is straightforward.

 \bip

 {\bf 3.7. Geometric interpretation of $\cP_\infty$-algebras}. A
  pre-Lie$^{\frak 2}_\infty$-structure on a dg vector space $V$ is a collection
  of linear maps,
  $$
  \mu_{k,p}: \odot^kV\ot \wedge^p V \rar V[1-p],\ \ k\geq 1,p\geq 0,
  $$
   satisfying a system
  of quadratic equations which can read off from equations 3.6.1 defining the cobar
  differential. There is a nice geometric way to describe these equations.

  \sip

  Let $\hat{V}$ be the smooth formal graded manifold isomorphic to the
  formal neighbourhood of $0$ in $V$, $T_{\hat{V}}$ the space of smooth vector fields,
 and
  $\Omega^\bullet_{\hat{V}}$ the de Rham algebra on $\hat{V}$.

 \bip

 If $\{e_\al, \al=1,2,\ldots\}$
is a homogeneous basis of $V$, then the associated dual basis $t^{\al}$,
 $|t^{\al}|=-|e_{\al}|$,
defines a coordinate system on $\hat{V}$. The collection of linear maps
$\{\mu_{k,0}\}_{k\geq 1}$ can be assembled into a germ, $\eth\in T_{\hat{V}}$, of
a degree 1 smooth vector field,
$$
 \eth:= \sum_{k=1}^{\infty}\frac{1}{k!} (-1)^{\epsilon}t^{\al_1}\cdots t^{\al_k}
\mu_{\al_1 \ldots \al_k}^{\ \ \ \ \ \ \be} \frac{\p}{\p t^{\be}}
$$
where
$$
\epsilon=\sum_{i=1}^k|e_{\al_i}|(1+\sum_{j=1}^i|e_{\al_j}|)
$$
and the numbers $\mu_{\al_1 \ldots \al_k}^{\ \ \ \ \ \ \be}$ are defined by
$$
\mu_{k,0}(e_{\al_1}, \ldots , e_{\al_k})=
\sum \mu_{\al_1 \ldots \al_k}^{\ \ \ \ \ \ \be}e_{\be}.
$$
We assume here and throughout the paper  summation over repeated small
Greek indices.

For fixed $p\geq 1$, the collection of linear maps,
$\{\mu_{k,p}\}_{k \geq 1}$ can be assembled into a germ,
$\Gamma_p\in T_{\hat{V}}\ot \Omega^p_{\hat{V}}$, $|\Gamma_p|=1$, of
a  smooth tangent vector valued differential $p$-form  on $\hat{V}$,
$$
 \Gamma_p:= \sum_{k=1}^{\infty}\frac{1}{k!p!} (-1)^{\epsilon}t^{\al_1}\cdots t^{\al_k}
\mu_{\al_1 \ldots \al_k, \be_1\ldots \be_p}^{\ \ \ \ \ \ \ \ \ \ \ \ \ \ \ \ga}
\frac{\p}{\p t^\ga}\ot dt^{\be_1}
\wedge \ldots \wedge  dt^{\be_p}
$$
where
$$
\epsilon=\sum_{i=1}^k|e_{\al_i}|(2-p+\sum_{j=1}^i|e_{\al_j}|)
+ \sum_{i=1}^k(|e_{\be_i}|+1)\sum_{j=i+1}^k|e_{\be_j}|
$$
and the numbers
$\mu_{\al_1 \ldots \al_k, \be_1\ldots \be_p}^{\ \ \ \ \ \ \ \ \ \ \ \ \ \ \ \ga}$
 are defined by
$$
\mu_{k,p}(e_{\al_1}, \ldots , e_{\al_k},
e_{\be_1}, \ldots , e_{\be_p})
=
\sum_\ga
\mu_{\al_1 \ldots \al_k, \be_1\ldots \be_p}^{\ \ \ \ \ \ \ \ \ \ \ \ \ \ \ \ga}
\frac{\p}{\p t^\ga}.
$$

\bip

{\bf 3.7.1. Proposition.} {\em A collection of linear maps,
$\{\mu_{k,p}: \odot^kV\ot \wedge^p V \rar V[1-p]\}_{k\geq 1,p\geq 0}$,
defines a structure of  pre-Lie$^{\frak 2}_\infty$-algebra
on $V$ if and only if the associated degree one smooth vector field $\eth$ and the
tangent vector  valued differential form,
$$
\Gamma := \sum_{p\geq 1} \Gamma_p\, \in T_{\hat{V}}\ot \Omega^\bullet_{\hat{V}}
$$
satisfy the equations,
$$
[\eth, \eth]=0
$$
and
$$
Lie_{\eth}\Gamma + \frac{1}{2}[\Gamma\bullet \Gamma]=0,
$$
where $Lie_\eth$ stands for the Lie derivative along the vector
field $\eth$ and $[\, \bullet \, ]$ for the  Fr\"olicher-Nijenhuis
brackets.}

\bip

In particular, if $V$ is finite dimensional and concentrated in degree zero, then
the only non-zero summand in $\Gamma$ is $\Gamma_1\in
T_{\hat{V}}\ot \Omega^1_{\hat{V}}$. According to 3.7.1, this endomorphism
of the tangent bundle $T_{\hat{V}}$ makes $V$ into a pre-Lie$^{\frak 2}_\infty$-algebra
if and only if the associated
 Nijenhuis tensor $N_{\Gamma_1}=[\Gamma_1\bullet \Gamma_1]$
vanishes.

\bip

{\bf 4. Further enlargement.} From the geometric viewpoint
there a small imperfection in the structure
of pre-Lie$^{\frak 2}_\infty$-algebras: the associated section $\Gamma$
of $T_{\hat{V}}\ot \Omega^\bullet_{\hat{V}}$ is required to be order $\geq 1$ in
the coordinate $t^\al$ so that constant sections are excluded.
In fact, the sections $\Gamma$ which are precisely linear in $t^\al$ together with
$\eth=0$
are in one-to-one correspondence to pre-Lie$_\infty$ structures on $V$, cf.\ \cite{CL}.

\bip

It is not a problem to fix that imperfection via a further enlargement
of the notion of pre-Lie algebra.

\bip

{\bf 4.1. Homologically trivial dg Lie algebras.} Such dg Lie algebras, $(\fg,d)$,
fit into an exact sequence,
$$
0\lon \ker d \lon \fg \stackrel{d}{\lon} \Img d \lon 0,
$$
or, denoting $V=\ker d[1]$, into the following one,
\Beqr\label{ext}
0\lon V[-1] \lon \fg \lon V\lon 0,
\Eeqr
which has all arrows of degree $0$.

\bip

{\bf 4.1.1. Definition I.} A {\em Nijenhuis algebra}, or shortly an $N$-{\em algebra},
  is a homologically trivial dg
Lie algebra, $(\fg, [\ , \, ], d)$, together with a fixed splitting,
$\fg=V\oplus V[-1]$, of the extension (\ref{ext}).

\sip

{\bf 4.1.2. Definition II.}  An {\em $N$-algebra} is a graded vector space $V$
together with
three binary operations,
$$
\Ba{rccc}
\circ: &  V\ot V & \lon & V \\
       & a\ot b    & \lon & a\circ b
\Ea
 , \
\Ba{rccc}
[\, \bullet\, ]: & \odot^2 V& \lon & V[1] \\
       & a\ot b    & \lon & (-1)^{|b|}[a\bullet b]
\Ea
 , \
\Ba{rccc}
{\star}: & \wedge^2 V& \lon & V[-1] \\
       & a\ot b    & \lon & (-1)^{|b|} a\star b
\Ea
$$
such that identities (1) and (2) of Definition 3.1 as well as the following ones,
$$
(a\circ b)\circ c -  (-1)^{|b||c|}(a\circ b)\circ c - a\circ (b\circ c)
+  (-1)^{|b||c|} a\circ (c\circ b)=\hspace{5cm}
$$
$$
\hspace{5cm}
 (-1)^{|b|} [a\bullet (b\star c) -(-1)^{|b|}[a\bullet b] \star c +
(-1)^{|b||c|+|c|}[a\bullet c]\star b,
$$

$$
(1+ \sigma + \sigma^2)\left\{ (a\circ b)\star c -(-1)^{|a||b|} (b\circ a)\star c +
(-1)^{|b|} (a\star b)\circ c \right\}=0,
$$
hold for any $a,b,c\in V$. Here $\sigma$ stands for the cyclic permutation of the
letters in the word $abc$.

\bip

{\bf 4.1.3. Proposition.} {\em Definitions I and II are equivalent}.
\bip

\Proof A binary operation,
$$
[\ ,\, ]: \wedge^2 (V\oplus V[-1])\lon V\oplus V[-1]
$$
satisfies the Leibniz rule with respect to the differential
$$
\Ba{rccc}
d: &  V\oplus V[-1] & \lon &  V\oplus V[-1]\\
       & a + \Pi b    & \lon & 0 +  \Pi a
\Ea
$$
if and only if it is of the form (cf.\ 3.1.2),
\Beqrn
[a,b]&:=& a\circ b - (-1)^{|a||b|} b\circ a + (-1)^{|a|}\Pi a\star b \\
\left[\Pi a, b\right] &:=& - (-1)^{|a|}[a\bullet b] +\Pi a\circ b\\
\left[\Pi a, \Pi b\right] &:=& \Pi\left[a\bullet b\right]
\Eeqrn
for some linear maps $\circ:  V\ot V  \rar V$,
$[\, \bullet\, ]: \odot^2 V \rar V[1]$ and
$\star: \wedge^2 V \rar V[-1]$. Then the Jacobi identities for $[\ ,\, ]$
get transformed into the identities of Definition II.
\hfill $\Box$

\bip

{\bf 4.1.4. SUSY transformations.} The set of all splittings of  extension
(\ref{ext}) is a principal homogeneous space over the Abelian group $V\ot V^*[-1]$.
Hence the latter acts as a kind of supersymmetry transformation on  $N$-algebra
structures: if $(V, [\, \bullet\, ], \circ,\star)$ is an $N$ algebra, then,
for any $f\in V\ot V^*[-1]$, the data $(V, [\, \bullet_f\, ], \circ_f,\star_f)$
is again an $N$-algebra, where
\Beqrn
[a\bullet_f b] &:=& [a\bullet b],\\
a\circ_f b &:=& a\circ b + [a\bullet f(b)] + (-1)^{|a|}f([a\bullet b],\\
a\star_f b &:=&
a\star b + (-1)^{|a|}\{f(a)\circ b - f(a\circ b)\}- (-1)^{(|a|+1)|b|}\{
f(b)\circ a - f(b\circ a)\}\\
&&
+ (-1)^{|a|}[f(a)\bullet f(b)] -  f[f(a)\bullet b]- (-1)^{|a|}f([a\bullet f(b)]),
\Eeqrn
for any $a,b,c\in V$.

\sip

The category of $N$-algebras with morphisms enlarged to include supersymmetry
transformations is equivalent to the category of homologically trivial
dg Lie algebras.

\bip

{\bf 4.2. Operad of $N$-algebras}. Let us denote this operad by $N$.
It is quadratic, $N={Free(\cE)}/{<R>}$, with
$\cE(2):=k[\bS_2][0]\oplus {\mathbf 1}_2[-1] \oplus sgn_2[1]$. Let the corollas
 $$
 
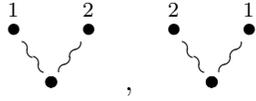
\begin{xy}
 <5mm,0cm>*{\bullet};<0cm,7mm>*{\bullet}**@{~},
 <5mm,0cm>*{\bullet};<10mm,7mm>*{\bullet}**@{~},
 <5mm,0cm>*{\bullet};<0mm,9mm>*{^1}**@{},
 <5mm,0cm>*{\bullet};<10mm,9mm>*{^2}**@{},
 \end{xy}
 \ \ \ , \ \ \
 \begin{xy}
 <5mm,0cm>*{\bullet};<0cm,7mm>*{\bullet}**@{~},
 <5mm,0cm>*{\bullet};<10mm,7mm>*{\bullet}**@{~},
 <5mm,0cm>*{\bullet};<0mm,9mm>*{^2}**@{},
 <5mm,0cm>*{\bullet};<10mm,9mm>*{^1}**@{},
 \end{xy} \ \  \ \
 $$
 stand for the basis of $k[\bS_2][0]$ and the corollas
$$
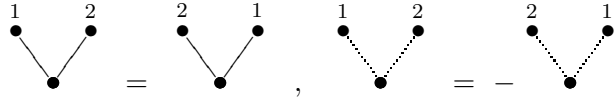
\begin{xy}
 <5mm,0cm>*{\bullet};<0cm,7mm>*{\bullet}**@{-},
 <5mm,0cm>*{\bullet};<10mm,7mm>*{\bullet}**@{-},
 <5mm,0cm>*{\bullet};<0mm,9mm>*{^1}**@{},
 <5mm,0cm>*{\bullet};<10mm,9mm>*{^2}**@{},
 \end{xy}
 \ \ =
\ \
\begin{xy}
 <5mm,0cm>*{\bullet};<0cm,7mm>*{\bullet}**@{-},
 <5mm,0cm>*{\bullet};<10mm,7mm>*{\bullet}**@{-},
 <5mm,0cm>*{\bullet};<0mm,9mm>*{^2}**@{},
 <5mm,0cm>*{\bullet};<10mm,9mm>*{^1}**@{},
 \end{xy}
\ \ \ ,
\ \ \
\begin{xy}
 <5mm,0cm>*{\bullet};<0cm,7mm>*{\bullet}**@{.},
 <5mm,0cm>*{\bullet};<10mm,7mm>*{\bullet}**@{.},
 <5mm,0cm>*{\bullet};<0mm,9mm>*{^1}**@{},
 <5mm,0cm>*{\bullet};<10mm,9mm>*{^2}**@{},
 \end{xy}
 \ \ =
\ - \
\begin{xy}
 <5mm,0cm>*{\bullet};<0cm,7mm>*{\bullet}**@{.},
 <5mm,0cm>*{\bullet};<10mm,7mm>*{\bullet}**@{.},
 <5mm,0cm>*{\bullet};<0mm,9mm>*{^2}**@{},
 <5mm,0cm>*{\bullet};<10mm,9mm>*{^1}**@{},
 \end{xy}
$$
for the basis vectors of  the ode-dimensional representations
$ {\mathbf 1}_2[-1]$ and $sgn_2[1]$ respectively.
Then the ideal ${<R>}$ is generated by  (\ref{A}), (\ref{B})
as well as
the following elements  of $Free(\cE)(3)$,
$$
 \begin{xy}
 <5mm,0mm>*{\bullet};<0cm,7mm>*{\bullet}**@{~},
 <5mm,0cm>*{\bullet};<10mm,7mm>*{\bullet}**@{~},
 <0mm,7mm>*{\bullet};<-5mm,14mm>*{\bullet}**@{~},
 <0mm,7mm>*{\bullet};<5mm,14mm>*{\bullet}**@{~},
 <5mm,0cm>*{\bullet};<-5mm,16.5mm>*{^{i_1}}**@{},
 <5mm,0cm>*{\bullet};<5mm,16.5mm>*{^{i_2}}**@{},
 <5mm,0cm>*{\bullet};<11mm,9.3mm>*{^{i_3}}**@{},
 \end{xy}
 -
  \begin{xy}
 <5mm,0mm>*{\bullet};<0cm,7mm>*{\bullet}**@{~},
 <5mm,0cm>*{\bullet};<10mm,7mm>*{\bullet}**@{~},
 <0mm,7mm>*{\bullet};<-5mm,14mm>*{\bullet}**@{~},
 <0mm,7mm>*{\bullet};<5mm,14mm>*{\bullet}**@{~},
 <5mm,0cm>*{\bullet};<-5mm,16.5mm>*{^{i_1}}**@{},
 <5mm,0cm>*{\bullet};<5mm,16.5mm>*{^{i_3}}**@{},
 <5mm,0cm>*{\bullet};<11mm,9.3mm>*{^{i_2}}**@{},
 \end{xy}
 -
 \begin{xy}
 <5mm,0mm>*{\bullet};<0cm,7mm>*{\bullet}**@{~},
 <5mm,0cm>*{\bullet};<10mm,7mm>*{\bullet}**@{~},
 <10mm,7mm>*{\bullet};<5mm,14mm>*{\bullet}**@{~},
 <10mm,7mm>*{\bullet};<15mm,14mm>*{\bullet}**@{~},
 <5mm,0cm>*{\bullet};<-1mm,9.3mm>*{^{i_1}}**@{},
 <5mm,0cm>*{\bullet};<5mm,16.5mm>*{^{i_2}}**@{},
 <5mm,0cm>*{\bullet};<16mm,16.5mm>*{^{i_3}}**@{},
 \end{xy}
+
\begin{xy}
 <5mm,0mm>*{\bullet};<0cm,7mm>*{\bullet}**@{~},
 <5mm,0cm>*{\bullet};<10mm,7mm>*{\bullet}**@{~},
 <10mm,7mm>*{\bullet};<5mm,14mm>*{\bullet}**@{~},
 <10mm,7mm>*{\bullet};<15mm,14mm>*{\bullet}**@{~},
 <5mm,0cm>*{\bullet};<-1mm,9.3mm>*{^{i_1}}**@{},
 <5mm,0cm>*{\bullet};<5mm,16.5mm>*{^{i_3}}**@{},
 <5mm,0cm>*{\bullet};<16mm,16.5mm>*{^{i_2}}**@{},
 \end{xy}
+
\begin{xy}
 <5mm,0mm>*{\bullet};<0cm,7mm>*{\bullet}**@{.},
 <5mm,0cm>*{\bullet};<10mm,7mm>*{\bullet}**@{.},
 <0mm,7mm>*{\bullet};<-5mm,14mm>*{\bullet}**@{-},
 <0mm,7mm>*{\bullet};<5mm,14mm>*{\bullet}**@{-},
 <5mm,0cm>*{\bullet};<-5mm,16.5mm>*{^{i_1}}**@{},
 <5mm,0cm>*{\bullet};<5mm,16.5mm>*{^{i_2}}**@{},
 <5mm,0cm>*{\bullet};<11mm,9.3mm>*{^{i_3}}**@{},
 \end{xy}
 -
 \begin{xy}
 <5mm,0mm>*{\bullet};<0cm,7mm>*{\bullet}**@{.},
 <5mm,0cm>*{\bullet};<10mm,7mm>*{\bullet}**@{.},
 <0mm,7mm>*{\bullet};<-5mm,14mm>*{\bullet}**@{-},
 <0mm,7mm>*{\bullet};<5mm,14mm>*{\bullet}**@{-},
 <5mm,0cm>*{\bullet};<-5mm,16.5mm>*{^{i_1}}**@{},
 <5mm,0cm>*{\bullet};<5mm,16.5mm>*{^{i_3}}**@{},
 <5mm,0cm>*{\bullet};<11mm,9.3mm>*{^{i_2}}**@{},
 \end{xy}
 -
 \begin{xy}
 <5mm,0mm>*{\bullet};<0cm,7mm>*{\bullet}**@{-},
 <5mm,0cm>*{\bullet};<10mm,7mm>*{\bullet}**@{-},
 <10mm,7mm>*{\bullet};<5mm,14mm>*{\bullet}**@{.},
 <10mm,7mm>*{\bullet};<15mm,14mm>*{\bullet}**@{.},
 <5mm,0cm>*{\bullet};<-1mm,9.3mm>*{^{i_1}}**@{},
 <5mm,0cm>*{\bullet};<5mm,16.5mm>*{^{i_2}}**@{},
 <5mm,0cm>*{\bullet};<16mm,16.5mm>*{^{i_3}}**@{},
 \end{xy}
 $$
$$
(1+\sigma+\sigma^2)(\
 \begin{xy}
 <5mm,0mm>*{\bullet};<0cm,7mm>*{\bullet}**@{~},
 <5mm,0cm>*{\bullet};<10mm,7mm>*{\bullet}**@{~},
 <0mm,7mm>*{\bullet};<-5mm,14mm>*{\bullet}**@{.},
 <0mm,7mm>*{\bullet};<5mm,14mm>*{\bullet}**@{.},
 <5mm,0cm>*{\bullet};<-5.4mm,16.3mm>*{^1}**@{},
 <5mm,0cm>*{\bullet};<5mm,16.3mm>*{^2}**@{},
 <5mm,0cm>*{\bullet};<10.6mm,9.1mm>*{^3}**@{},
 \end{xy}\ \ + \ \
\begin{xy}
 <5mm,0mm>*{\bullet};<0cm,7mm>*{\bullet}**@{.},
 <5mm,0cm>*{\bullet};<10mm,7mm>*{\bullet}**@{.},
 <0mm,7mm>*{\bullet};<-5mm,14mm>*{\bullet}**@{~},
 <0mm,7mm>*{\bullet};<5mm,14mm>*{\bullet}**@{~},
 <5mm,0cm>*{\bullet};<-5.4mm,16.3mm>*{^1}**@{},
 <5mm,0cm>*{\bullet};<5mm,16.3mm>*{^2}**@{},
 <5mm,0cm>*{\bullet};<10.6mm,9.1mm>*{^3}**@{},
 \end{xy}\ \ - \ \
 \begin{xy}
 <5mm,0mm>*{\bullet};<0cm,7mm>*{\bullet}**@{.},
 <5mm,0cm>*{\bullet};<10mm,7mm>*{\bullet}**@{.},
 <0mm,7mm>*{\bullet};<-5mm,14mm>*{\bullet}**@{~},
 <0mm,7mm>*{\bullet};<5mm,14mm>*{\bullet}**@{~},
 <5mm,0cm>*{\bullet};<-5.4mm,16.3mm>*{^2}**@{},
 <5mm,0cm>*{\bullet};<5mm,16.3mm>*{^1}**@{},
 <5mm,0cm>*{\bullet};<10.6mm,9.1mm>*{^3}**@{},
 \end{xy}
\ \ ) ,
$$
where $\sigma$ is the cyclic permutation of $(1,2,3)$.

 \bip

{\bf 4.3. Koszul dual of $N$.} By definition,
$N^!={Free(\check{\cE)}}/{<R^\bot>}$, with the only non-vanishing bit
of the $\bS$-module $\check{\cE}$ being
$\check{\cE}(2):=k[\bS_2][0]  \oplus sgn_2[1]\oplus {\mathbf 1}_2[-1]$. Let the corollas
 $$
 
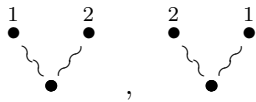
\begin{xy}
 <5mm,0cm>*{\bullet};<0cm,7mm>*{\bullet}**@{~},
 <5mm,0cm>*{\bullet};<10mm,7mm>*{\bullet}**@{~},
 <5mm,0cm>*{\bullet};<0mm,9mm>*{^1}**@{},
 <5mm,0cm>*{\bullet};<10mm,9mm>*{^2}**@{},
 \end{xy}
 \ \ \ , \ \ \
 \begin{xy}
 <5mm,0cm>*{\bullet};<0cm,7mm>*{\bullet}**@{~},
 <5mm,0cm>*{\bullet};<10mm,7mm>*{\bullet}**@{~},
 <5mm,0cm>*{\bullet};<0mm,9mm>*{^2}**@{},
 <5mm,0cm>*{\bullet};<10mm,9mm>*{^1}**@{},
 \end{xy} \ \  \ \
 $$
 stand for the basis of $k[\bS_2][0]$ while the corollas
$$
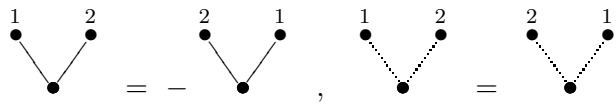
\begin{xy}
 <5mm,0cm>*{\bullet};<0cm,7mm>*{\bullet}**@{-},
 <5mm,0cm>*{\bullet};<10mm,7mm>*{\bullet}**@{-},
 <5mm,0cm>*{\bullet};<0mm,9mm>*{^1}**@{},
 <5mm,0cm>*{\bullet};<10mm,9mm>*{^2}**@{},
 \end{xy}
 \ \ =
\ - \
\begin{xy}
 <5mm,0cm>*{\bullet};<0cm,7mm>*{\bullet}**@{-},
 <5mm,0cm>*{\bullet};<10mm,7mm>*{\bullet}**@{-},
 <5mm,0cm>*{\bullet};<0mm,9mm>*{^2}**@{},
 <5mm,0cm>*{\bullet};<10mm,9mm>*{^1}**@{},
 \end{xy}
\ \ \ ,
\ \ \
\begin{xy}
 <5mm,0cm>*{\bullet};<0cm,7mm>*{\bullet}**@{.},
 <5mm,0cm>*{\bullet};<10mm,7mm>*{\bullet}**@{.},
 <5mm,0cm>*{\bullet};<0mm,9mm>*{^1}**@{},
 <5mm,0cm>*{\bullet};<10mm,9mm>*{^2}**@{},
 \end{xy}
 \ \ =
\  \
\begin{xy}
 <5mm,0cm>*{\bullet};<0cm,7mm>*{\bullet}**@{.},
 <5mm,0cm>*{\bullet};<10mm,7mm>*{\bullet}**@{.},
 <5mm,0cm>*{\bullet};<0mm,9mm>*{^2}**@{},
 <5mm,0cm>*{\bullet};<10mm,9mm>*{^1}**@{},
 \end{xy}
$$
represent the basis vectors of $sgn_2[1]$ and
$ {\mathbf 1}_2[-1]$  respectively.
Then the ideal ${<R^\bot>}$ is generated by (\ref{C}), (\ref{D}), (\ref{E}) and
 following elements of $Free(\check{\cE})(3)$,
$$
 
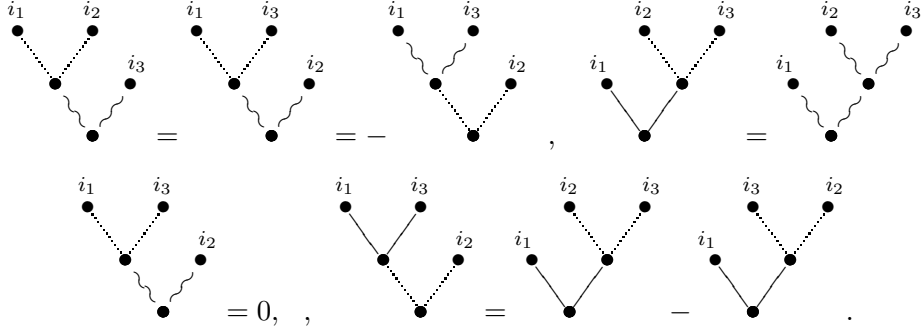
\begin{xy}
 <5mm,0mm>*{\bullet};<0cm,7mm>*{\bullet}**@{~},
 <5mm,0cm>*{\bullet};<10mm,7mm>*{\bullet}**@{~},
 <0mm,7mm>*{\bullet};<-5mm,14mm>*{\bullet}**@{.},
 <0mm,7mm>*{\bullet};<5mm,14mm>*{\bullet}**@{.},
 <5mm,0cm>*{\bullet};<-5mm,16.5mm>*{^{i_1}}**@{},
 <5mm,0cm>*{\bullet};<5mm,16.5mm>*{^{i_2}}**@{},
 <5mm,0cm>*{\bullet};<11mm,9.3mm>*{^{i_3}}**@{},
 \end{xy}
 =
\begin{xy}
 <5mm,0mm>*{\bullet};<0cm,7mm>*{\bullet}**@{~},
 <5mm,0cm>*{\bullet};<10mm,7mm>*{\bullet}**@{~},
 <0mm,7mm>*{\bullet};<-5mm,14mm>*{\bullet}**@{.},
 <0mm,7mm>*{\bullet};<5mm,14mm>*{\bullet}**@{.},
 <5mm,0cm>*{\bullet};<-5mm,16.5mm>*{^{i_1}}**@{},
 <5mm,0cm>*{\bullet};<5mm,16.5mm>*{^{i_3}}**@{},
 <5mm,0cm>*{\bullet};<11mm,9.3mm>*{^{i_2}}**@{},
 \end{xy}
=
-
\begin{xy}
 <5mm,0mm>*{\bullet};<0cm,7mm>*{\bullet}**@{.},
 <5mm,0cm>*{\bullet};<10mm,7mm>*{\bullet}**@{.},
 <0mm,7mm>*{\bullet};<-5mm,14mm>*{\bullet}**@{~},
 <0mm,7mm>*{\bullet};<5mm,14mm>*{\bullet}**@{~},
 <5mm,0cm>*{\bullet};<-5mm,16.5mm>*{^{i_1}}**@{},
 <5mm,0cm>*{\bullet};<5mm,16.5mm>*{^{i_3}}**@{},
 <5mm,0cm>*{\bullet};<11mm,9.3mm>*{^{i_2}}**@{},
 \end{xy}
\ \ , \ \ \
\begin{xy}
 <5mm,0mm>*{\bullet};<0cm,7mm>*{\bullet}**@{-},
 <5mm,0cm>*{\bullet};<10mm,7mm>*{\bullet}**@{-},
 <10mm,7mm>*{\bullet};<5mm,14mm>*{\bullet}**@{.},
 <10mm,7mm>*{\bullet};<15mm,14mm>*{\bullet}**@{.},
 <5mm,0cm>*{\bullet};<-1mm,9.3mm>*{^{i_1}}**@{},
 <5mm,0cm>*{\bullet};<5mm,16.5mm>*{^{i_2}}**@{},
 <5mm,0cm>*{\bullet};<16mm,16.5mm>*{^{i_3}}**@{},
 \end{xy}
 =
 \begin{xy}
 <5mm,0mm>*{\bullet};<0cm,7mm>*{\bullet}**@{~},
 <5mm,0cm>*{\bullet};<10mm,7mm>*{\bullet}**@{~},
 <10mm,7mm>*{\bullet};<5mm,14mm>*{\bullet}**@{~},
 <10mm,7mm>*{\bullet};<15mm,14mm>*{\bullet}**@{~},
 <5mm,0cm>*{\bullet};<-1mm,9.3mm>*{^{i_1}}**@{},
 <5mm,0cm>*{\bullet};<5mm,16.5mm>*{^{i_2}}**@{},
 <5mm,0cm>*{\bullet};<16mm,16.5mm>*{^{i_3}}**@{},
 \end{xy}
$$
$$
\begin{xy}
 <5mm,0mm>*{\bullet};<0cm,7mm>*{\bullet}**@{~},
 <5mm,0cm>*{\bullet};<10mm,7mm>*{\bullet}**@{~},
 <0mm,7mm>*{\bullet};<-5mm,14mm>*{\bullet}**@{.},
 <0mm,7mm>*{\bullet};<5mm,14mm>*{\bullet}**@{.},
 <5mm,0cm>*{\bullet};<-5mm,16.5mm>*{^{i_1}}**@{},
 <5mm,0cm>*{\bullet};<5mm,16.5mm>*{^{i_3}}**@{},
 <5mm,0cm>*{\bullet};<11mm,9.3mm>*{^{i_2}}**@{},
 \end{xy}
=0,
\ \ , \ \
\begin{xy}
 <5mm,0mm>*{\bullet};<0cm,7mm>*{\bullet}**@{.},
 <5mm,0cm>*{\bullet};<10mm,7mm>*{\bullet}**@{.},
 <0mm,7mm>*{\bullet};<-5mm,14mm>*{\bullet}**@{-},
 <0mm,7mm>*{\bullet};<5mm,14mm>*{\bullet}**@{-},
 <5mm,0cm>*{\bullet};<-5mm,16.5mm>*{^{i_1}}**@{},
 <5mm,0cm>*{\bullet};<5mm,16.5mm>*{^{i_3}}**@{},
 <5mm,0cm>*{\bullet};<11mm,9.3mm>*{^{i_2}}**@{},
 \end{xy}
 =
\begin{xy}
 <5mm,0mm>*{\bullet};<0cm,7mm>*{\bullet}**@{-},
 <5mm,0cm>*{\bullet};<10mm,7mm>*{\bullet}**@{-},
 <10mm,7mm>*{\bullet};<5mm,14mm>*{\bullet}**@{.},
 <10mm,7mm>*{\bullet};<15mm,14mm>*{\bullet}**@{.},
 <5mm,0cm>*{\bullet};<-1mm,9.3mm>*{^{i_1}}**@{},
 <5mm,0cm>*{\bullet};<5mm,16.5mm>*{^{i_2}}**@{},
 <5mm,0cm>*{\bullet};<16mm,16.5mm>*{^{i_3}}**@{},
 \end{xy}
 -
 \begin{xy}
 <5mm,0mm>*{\bullet};<0cm,7mm>*{\bullet}**@{-},
 <5mm,0cm>*{\bullet};<10mm,7mm>*{\bullet}**@{-},
 <10mm,7mm>*{\bullet};<5mm,14mm>*{\bullet}**@{.},
 <10mm,7mm>*{\bullet};<15mm,14mm>*{\bullet}**@{.},
 <5mm,0cm>*{\bullet};<-1mm,9.3mm>*{^{i_1}}**@{},
 <5mm,0cm>*{\bullet};<5mm,16.5mm>*{^{i_3}}**@{},
 <5mm,0cm>*{\bullet};<16mm,16.5mm>*{^{i_2}}**@{},
 \end{xy}.
$$

\bip

{\bf 4.3.1. Example}. If $(A, \cdot, d)$ is a
differential graded commutative associative algebra, then $A[-1]$
has a natural structure of $N^!$-algebra with the products
(cf.\ Example 3.4.1)
\Beqrn
\Pi \al \circ \Pi\be & = & \Pi (\al\cdot d\be),\\
\Pi \al \bullet \Pi\be & = & \Pi (\al\cdot \be),\\
\Pi \al \star \Pi\be & = & (-1)^{|a|}\Pi (d\al\cdot d\be),
\Eeqrn
where $\circ$, $\bullet$ and $\star$ correspond, respectively,
to the generators
$$
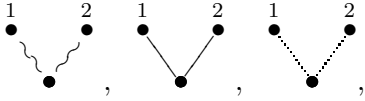
\begin{xy}
 <5mm,0cm>*{\bullet};<0cm,7mm>*{\bullet}**@{~},
 <5mm,0cm>*{\bullet};<10mm,7mm>*{\bullet}**@{~},
 <5mm,0cm>*{\bullet};<0mm,9mm>*{^1}**@{},
 <5mm,0cm>*{\bullet};<10mm,9mm>*{^2}**@{},
 \end{xy}
 \ ,
 \ \
 \begin{xy}
 <5mm,0cm>*{\bullet};<0cm,7mm>*{\bullet}**@{-},
 <5mm,0cm>*{\bullet};<10mm,7mm>*{\bullet}**@{-},
 <5mm,0cm>*{\bullet};<0mm,9mm>*{^1}**@{},
 <5mm,0cm>*{\bullet};<10mm,9mm>*{^2}**@{},
 \end{xy}
\ ,
\ \
\begin{xy}
 <5mm,0cm>*{\bullet};<0cm,7mm>*{\bullet}**@{.},
 <5mm,0cm>*{\bullet};<10mm,7mm>*{\bullet}**@{.},
 <5mm,0cm>*{\bullet};<0mm,9mm>*{^1}**@{},
 <5mm,0cm>*{\bullet};<10mm,9mm>*{^2}**@{},
 \end{xy},
$$
of the operad $N^!$.

\bip

{\bf 5. Geometric description of Nijenhuis$_\infty$ algebras}. The cobar construction,
$N_\infty:={\mathbf D}J^!$ is the free operad on the
$\bS$-module $\bigoplus_{p=0}^{n}
{\mathrm Ind}^{\bS_n}_{\bS_{n-p}\times \bS_{p}}
{\mathbf 1}_{n-p}\ot sgn_{p}[p-1]$. The difference from the pre-Lie$^{\frak 2}$
case in Section 3.6 is an extra summand corresponding to $p=n$. The cobar differential
has the form of the Proposition 3.6.1 with the summation ranges  extended to
include
\Bi
\item
the case $\# I_2=0$ in the first sum and
\item
the case  $\# I_1=0$ in the second sum.
\Ei
Hence an $N_\infty$-structure on a graded vector space $V$ can be described as
\Bi
\item a collection of linear maps,
$\{\mu_{k,p}: \odot^k V\ot \wedge^p V \rar V[1-p]\}_{k\geq 1, p\geq 0}$,
exactly the same as in the case of pre-Lie$^{\frak 2}$ algebras, plus
\item  an additional collection of linear maps
$\{\mu_{0,p}: \wedge^p V \rar V[1-p]\}_{p\geq 1}\}$,
\Ei
satisfying the equations which can be interpreted as in
Proposition~3.7.1: the differential form
$$
\Gamma := \sum_{p\geq 1}\sum_{k=0}^{\infty}\frac{1}{k!p!} (-1)^{\epsilon}
t^{\al_1}\cdots t^{\al_k}
\mu_{\al_1 \ldots \al_k, \be_1\ldots \be_p}^{\ \ \ \ \ \ \ \ \ \ \ \ \ \ \ \ga}
\frac{\p}{\p t^\ga}\ot dt^{\be_1}
\wedge \ldots \wedge  dt^{\be_p}\in T_{\hat{V}}\ot \Omega^\bullet_{\hat{V}}
$$
satisfies $
Lie_{\eth}\Gamma + \frac{1}{2}[\Gamma\bullet \Gamma]=0$. Representations of
$N_\infty$ in vector spaces concentrated in degree 0 are obviously
in one-to-one correspondence
with formal germs of Nijenhuis geometric structures which vanish at $0$.

\sip

There is, however, another more coherent picture of $N_\infty$ algebras.

\bip

{\bf 5.1. Reminder on the homotopy classification of $L_\infty$ algebras \cite{K}}.
Let $L$ be an operad of Lie algebras and $L_\infty$ its cobar construction (which
happens to be the minimal resolution of $L$ as the latter is Koszul).
There is a one-to-one correspondence between
$L_\infty$ algebra structures on a dg vector space $V$
and degree one vector fields $\eth$ on $\hat{V}$, the formal neighbourhood of zero in $V$,
which vanish at $0$ and satisfy
the integrability condition $[\eth,\eth]=0$. The data $(\hat{V},\eth)$ is often called
a smooth formal dg manifold.

\sip

Hence the homotopy theory of $L_\infty$ algebras is the same as the homotopy
theory of formal dg manifolds.

\sip

A $L_\infty$ algebra $(\hat{V}, \eth)$, is called {\em minimal}\, if the first, $\eth_1$,
Taylor coefficient of the homological vector field $\eth$ vanishes. It is called
{\em linear contractible}\, if the higher Taylor coefficients $\eth_{\geq 2}$ vanish
and the first one $\eth_1$ has trivial cohomology when viewed as a differential in $V$.

\sip

According to Kontsevich \cite{K}, any dg manifold is $L_\infty$-isomorphic to
the direct product of a minimal and of a linear contractible one. A dg manifold
is called {\em contractible}\, if it is $L_\infty$-isomorphic to
 a linear contractible one.

\bip

{\bf 5.2. Nijenhuis$_\infty$ versus contractible dg manifolds}.
Let $\cM$ be a formal neighbourhood  of zero in $V\oplus V[-1]$. The standard differential $d$
on $V\oplus V[-1]$ induces a homological vector field on $\cM$ which we denote by the same
letter. Thus $(\cM,d)$ is a linear contractible dg manifold.

\sip
Let $i:  T_{\hat{V}}\ot \Omega^\bullet_{\hat{V}}\rar T_\cM$ be degree $-1$ linear
map of subsection 3.1.4 which sends sections of  $T_{\hat{V}}\ot \Omega^\bullet_{\hat{V}}$
into vertical vector fields on the fibration $\cM\rar \hat{V}$.
Consider the degree zero map
$$
\Ba{rccc}
\Psi: &   T_{\hat{V}}\ot \Omega^\bullet_{\hat{V}}& \lon & T_\cM \\
&  X &\lon & [d,i(X)].
\Ea
$$
The equation  $
Lie_{\eth}\Gamma + \frac{1}{2}[\Gamma\bullet \Gamma]=0$ implies that
 $\hat{\eth}=\Psi(\eth + \Gamma)$ is  a homological vector field
on $\cM$ which commutes with $d$.
Thus we obtain the following statement.

\bip

{\bf 5.2.1. Theorem}. {\em There is a one-to-one correspondence between
 Nijenhuis$_\infty$ structures in a dg vector space $V$ and homological vector fields
 on $\cM$ which commute with $d$.}

\bip

A  Nijenhuis$_\infty$ algebra is called {\em minimal}\, if
the homological vector
field $\eth$ (or, equivalently, $\hat{\eth}$)
 is minimal, i.e.\ its first Taylor coefficient vanishes.
 A more invariant way
to formulate this condition is to say that $\hat{\eth} I\subset I^2$ where
$I$ is the ideal of the distinguished  point $0\in \cM$.
Then it follows from Theorem~5.2.1 that
 $d + \hat{\eth}$ is a homological vector field on $\cM$ whose first
Taylor coefficient coincides with $d$.

\bip

{\bf 5.2.2. Corollary}. {\em There is a one-to-one correspondence}
$$
\left\{\mbox{minimal Nijenhuis$_\infty$\ algebras}\right\}
\stackrel{1:1}\longleftrightarrow
\left\{\mbox{contractible dg manifolds}\right\}.
$$

\bip

As the operad $N_\infty$ is cofibrant,  any Nijenhuis$_\infty$ algebra
is homotopy equivalent to a minimal one. Hence

\bip

{\bf 5.2.3. Corollary}. {\em Any
Nijenhuis$_\infty$\ algebra is homotopy equivalent to a
contractible dg manifold}.

\bip

\pagebreak

 \bip

  {\small

  \end{document}